\documentclass[12pt]{article} 

\usepackage{latexsym, 
amscd, 
graphicx, amssymb}

\newtheorem{thm}{Theorem}[section]

\newtheorem{prop}[thm]{Proposition}
\newtheorem{cor}[thm]{Corollary}
\newtheorem{lemma}[thm]{Lemma}
\newtheorem{rema}[thm]{Remark}

\newcommand{\halmos}{\rule{1ex}{1.4ex}}

\newcommand{\text}[1]{\mbox{\rm #1}}

\newcommand{\nn}{\nonumber \\}

 \newcommand{\res}{\mbox{\rm Res}}
 
\renewcommand{\hom}{\mbox{\rm Hom}}

 \newcommand{\pf}{{\it Proof.}\hspace{2ex}}
 
 \newcommand{\epfv}{\hspace*{\fill}\mbox{$\halmos$}\vspace{1em}}

\newcommand{\tr}{\mbox{\rm Tr}}

\newcommand{\A}{\mathcal{A}}
\newcommand{\Y}{\mathcal{Y}}
\newcommand{\C}{\mathbb{C}}
\newcommand{\Z}{\mathbb{Z}}
\newcommand{\R}{\mathbb{R}}

\newcommand{\N}{\mathbb{N}}

\newcommand{\V}{\mathcal{V}}

\makeatletter
\newlength{\@pxlwd} \newlength{\@rulewd} \newlength{\@pxlht}
\catcode`.=\active \catcode`B=\active \catcode`:=\active \catcode`|=\active
\def\sprite#1(#2,#3)[#4,#5]{
   \edef\@sprbox{\expandafter\@cdr\string#1\@nil @box}
   \expandafter\newsavebox\csname\@sprbox\endcsname
   \edef#1{\expandafter\usebox\csname\@sprbox\endcsname}
   \expandafter\setbox\csname\@sprbox\endcsname =\hbox\bgroup
   \vbox\bgroup
  \catcode`.=\active\catcode`B=\active\catcode`:=\active\catcode`|=\active
      \@pxlwd=#4 \divide\@pxlwd by #3 \@rulewd=\@pxlwd
      \@pxlht=#5 \divide\@pxlht by #2
      \def .{\hskip \@pxlwd \ignorespaces}
      \def B{\@ifnextchar B{\advance\@rulewd by \@pxlwd}{\vrule
         height \@pxlht width \@rulewd depth 0 pt \@rulewd=\@pxlwd}}
      \def :{\hbox\bgroup\vrule height \@pxlht width 0pt depth
0pt\ignorespaces}
      \def |{\vrule height \@pxlht width 0pt depth 0pt\egroup
         \prevdepth= -1000 pt}
   }
\def\endsprite{\egroup\egroup}
\catcode`.=12 \catcode`B=11 \catcode`:=12 \catcode`|=12\relax
\makeatother

\def\hboxtr{\FormOfHboxtr} 
\sprite{\FormOfHboxtr}(25,25)[0.5 em, 1.2 ex] 
:BBBBBBBBBBBBBBBBBBBBBBBBB |
:BB......................B |
:B.B.....................B |
:B..B....................B |
:B...B...................B |
:B....B..................B |
:B.....B.................B |
:B......B................B |
:B.......B...............B |
:B........B..............B |
:B.........B.............B |
:B..........B............B |
:B...........B...........B |
:B............B..........B |
:B.............B.........B |
:B..............B........B |
:B...............B.......B |
:B................B......B |
:B.................B.....B |
:B..................B....B |
:B...................B...B |
:B....................B..B |
:B.....................B.B |
:B......................BB |
:BBBBBBBBBBBBBBBBBBBBBBBBB |
\endsprite

\sprite{\FormOfShboxtr}(25,25)[0.3 em, 0.72 ex] 
:BBBBBBBBBBBBBBBBBBBBBBBBB |
:BB......................B |
:B.B.....................B |
:B..B....................B |
:B...B...................B |
:B....B..................B |
:B.....B.................B |
:B......B................B |
:B.......B...............B |
:B........B..............B |
:B.........B.............B |
:B..........B............B |
:B...........B...........B |
:B............B..........B |
:B.............B.........B |
:B..............B........B |
:B...............B.......B |
:B................B......B |
:B.................B.....B |
:B..................B....B |
:B...................B...B |
:B....................B..B |
:B.....................B.B |
:B......................BB |
:BBBBBBBBBBBBBBBBBBBBBBBBB |
\endsprite

\title{ {\bf Rigidity and modularity of vertex tensor categories} }
\date{}
\author{Yi-Zhi Huang}

\begin{document}

\bibliographystyle{alpha}
\maketitle

\begin{abstract} 
Let $V$ be a simple vertex operator algebra satisfying the 
following conditions: (i) $V_{(n)}=0$ for $n<0$, 
$V_{(0)}=\mathbb{C}\mathbf{1}$ and 
$V'$ is isomorphic to $V$ as a $V$-module. (ii)
Every $\mathbb{N}$-gradable weak $V$-module is completely 
reducible. (iii) $V$ is $C_{2}$-cofinite. (In the presence 
of Condition (i), Conditions (ii) and (iii) are equivalent to
a single condition, namely, that every weak $V$-module is completely 
reducible.) Using the results obtained by the author 
in the formulation and proof of the general version of the Verlinde 
conjecture and in the proof of the Verlinde formula, we prove that the 
braided tensor category structure on the category
of $V$-modules is rigid, balanced and nondegenerate.
In particular, the category of $V$-modules 
has a natural structure of modular tensor category. 
We also prove that the tensor-categorical 
dimension of an irreducible $V$-module 
is the reciprocal of a suitable matrix element of the 
fusing isomorphism under a suitable basis. 
\end{abstract}

\renewcommand{\theequation}{\thesection.\arabic{equation}}
\renewcommand{\thethm}{\thesection.\arabic{thm}}
\setcounter{equation}{0}
\setcounter{thm}{0}
\setcounter{section}{-1}

\section{Introduction}

In the present paper, we prove the rigidity and modularity of 
the braided tensor category of modules for a vertex operator
algebra satisfying certain natural conditions (see below). 
Finding proofs of these properties has been an 
open problem for many years. Our proofs in this paper are based 
on the results obatined by the author \cite{H7}
in the formulation
and proof of the general version of the Verlinde 
conjecture and in the proof of the Verlinde formula.

In 1988, Moore and Seiberg \cite{MS1} \cite{MS2} derived a system of
polynomial equations from the axioms for rational conformal field
theories. They showed that the Verlinde conjecture \cite{V} is a
consequence of these equations. Inspired by an observation of Witten
on an analogy with Mac Lane's coherence, Moore and Seiberg \cite{MS2}
also demonstrated that the theory of these polynomial equations is
actually a conformal-field-theoretic analogue of the theory of the tensor
categories.  This work of Moore and Seiberg greatly advanced our
understanding of the structure of conformal field theories
and the name ``modular tensor category'' was suggested by
I. Frenkel for the
theory of these Moore-Seiberg equations.  Later, the precise notion of
modular tensor category was introduced to summarize the properties of
these polynomial equations and has played a central role in
the development of conformal field theories and three-dimensional
topological field theories.  See for example \cite{T} and \cite{BK} for
the theory of modular tensor categories, their applications and
references to many important works of mathematicians and
physicists.

Mathematically, Kazhdan and Lusztig \cite{KL1}--\cite{KL5} first
constructed a rigid braided tensor category structure on a suitable
(nonsemisimple) category of modules of a negative level for
an affine Lie algebra. Finkelberg \cite{F1} \cite{F2} transported these
braided tensor category structures to the category of integrable
highest weight modules of positive integral levels (with a few
exceptions) for the same affine Lie algebra. Direct constructions 
of these braided tensor category structures were also given by Lepowsky 
and the author 
\cite{HL5} based on the results in \cite{HL1}--\cite{HL4} \cite{H1}
and by Bakalov and Kirillov \cite{BK}. In the general case, for a
vertex operator algebra $V$ satisfying suitable conditions (weaker than
the conditions in the present paper), the braided tensor category
structure on the category of $V$-modules was constructed by Lepowsky and
the author and by the author in a series of papers
\cite{HL1}--\cite{HL4} \cite{H1} \cite{H5}.  This construction has been
generalized to the nonsemisimple (logarithmic) case by Lepowsky, Zhang
and the author recently \cite{HLZ1} \cite{HLZ2}, and as an application, 
a different construction of 
Kazhdan-Lusztig's braided tensor category structure has been given 
using this logarithmic theory by Zhang \cite{Zh} \cite{Zh2}. 

To prove that a semisimple braided tensor category actually carries
a modular tensor category structure, we need to prove that it is
rigid, balanced and nondegenerate. 
The balancing isomorphisms or twists in these braided tensor categories are
actually trivial to construct and the balancing axioms are easy to
prove.  On the other hand, the rigidity 
has been an open problem for many years
after the braided tensor category structure on the
category of modules for a vertex operator algebra 
satisfying the conditions in \cite{HL1}--\cite{HL4} \cite{H1}
was constructed.  
The main difficulty
is that from the construction, it is not clear why the numbers
determining the module maps given by the sequences in the axioms for the 
rigidity are not $0$. The nondegeneracy of the semisimple 
braided tensor category of
modules for a suitable general vertex operator algebra has also
been open. Even in the special case of  
the category of integrable
highest weight modules of positive integral levels for 
affine Lie algebras, 
as far as the author knows, there have been no proofs of 
the rigidity or the nondegeneracy property in the literature. 

In the present paper, we solve all these problems.  Let $V$ be a
simple vertex operator algebra satisfying the following conditions: 
(i) $V_{(n)}=0$ for $n<0$, 
$V_{(0)}=\mathbb{C}\mathbf{1}$ and 
$V'$ is isomorphic to $V$ as a $V$-module. (ii) Every
$\mathbb{N}$-gradable weak $V$-module is completely reducible. (iii) $V$
is $C_{2}$-cofinite. In an early version of the present paper and in 
the results announced in \cite{H8} and \cite{H9}, the proofs of the 
rigidity and the nondegeneracy property use a condition (Condition (i)
in the early version of the present paper)
slightly stronger than Condition (i) above. 
It requires that 
$V_{(n)}=0$ for $n<0$, $V_{(0)}=\mathbb{C}\mathbf{1}$ and 
for any irreducible $V$-module $W$, $W_{(0)}=0$.
But actually both proofs still work under Condition (i) 
in the present version\footnote{I am grateful to Liang Kong 
for pointing out
that the proof of the rigidity in the early version 
still works with Condition (i) in the present version.}. 
Also, by results of Li \cite{L} and 
Abe-Buhl-Dong \cite{ABD}, 
Conditions (ii) and (iii) are equivalent to a 
single condition that every weak $V$-module is completely 
reducible. In the present paper, using a consequence of 
a Verlinde formula
proved recently by the author in \cite{H7}, we show that the braided
tensor category structure on the category of $V$-modules is rigid. Using
a formula also obtained easily from this Verlinde formula (see \cite{H7}), we
prove that the semisimple rigid balanced braided tensor category
structure on the category of $V$-modules is nondegenerate.  In
particular, the category of $V$-modules has a natural structure of
modular tensor category. 

The results of the present paper have been announced 
in \cite{H8} and \cite{H9}. See also \cite{Le} for an exposition.

The paper is organized as follows: Section 1 is a review of the 
tensor product theory developed by Lepowsky and the author 
in \cite{HL1}--\cite{HL4}, \cite{HL6}, \cite{H1} and \cite{H5}. 
Section 2 is a review of the fusing and braiding matrices, the 
Verlinde conjecture and its consequences studied and proved by the 
author in \cite{H1}--\cite{H7}. See also \cite{H8} and \cite{H9} for 
expositions. In Section 3, we prove the rigidity of the braided tensor 
category of $V$-modules for a vertex operator algebra satisfying the 
three conditions mentioned above. We remark that 
we can introduce a notion of rigidity of vertex tensor categories and 
we actually have proved that the vertex tensor category structure 
on the category of $V$-modules is rigid in this sense. At the end of 
this section, we calculate the tensor-categorical 
dimension of irreducible $V$-modules explicitly. 
In Section 4, we first show that the semi-simple 
rigid braided tensor category
of $V$-modules with the obvious twists is a ribbon category.
Then we show that this ribbon category is 
nondegenerate. In particular, the category 
of $V$-modules is a modular tensor category. We also 
remark that 
we can introduce a notion of modular vertex tensor categories and 
we actually have proved that the category of $V$-modules is  
a modular vertex tensor category in this sense.

\paragraph{Acknowledgment} This research is partially supported 
by NSF grants DMS-0070800 and DMS-0401302. I am
grateful to the referee, to J. Lepowsky and especially to
L. Kong for their comments.

\renewcommand{\theequation}{\thesection.\arabic{equation}}
\renewcommand{\thethm}{\thesection.\arabic{thm}}
\setcounter{equation}{0}
\setcounter{thm}{0}

\section{Review of the tensor product theory}

In this section, we review the tensor product theory for modules
for a vertex operator algebra developed by Lepowsky and the author 
in \cite{HL1}--\cite{HL4}, \cite{HL6}, \cite{H1} and \cite{H5}.

Let $V$ be a simple vertex operator algebra and $C_{2}(V)$ the subspace of 
$V$ spanned by $u_{-2}v$ for $u, v\in V$. In the present paper,
we shall always assume that $V$ satisfies the 
following conditions:

\begin{enumerate}

\item $V_{(n)}=0$ for $n<0$, 
$V_{(0)}=\mathbb{C}\mathbf{1}$ and 
$V'$ is isomorphic to $V$ as a $V$-module.

\item Every $\mathbb{N}$-gradable weak $V$-module is completely 
reducible.

\item $V$ is $C_{2}$-cofinite, that is, $\dim V/C_{2}(V)<\infty$. 

\end{enumerate}

These conditions are all natural.  Condition 1 says that the vacuum is
unique in $V$ and the contragredient of $V$ does not give a
new irreducible $V$-module. As we mentioned in the 
introduction, Condition 1 is a consequence of the 
following stronger version of the uniqueness-of-the vacuum 
condition:  $V_{(n)}=0$ for $n<0$, $V_{(0)}=\mathbb{C}\mathbf{1}$ and 
for any irreducible $V$-module $W$, $W_{(0)}=0$. 
Note that finitely generated $\N$-gradable weak
$V$-modules are what naturally appear in the proofs of the theorems on
genus-zero and genus-one correlation functions. Thus Condition 2 is
natural and necessary because the Verlinde conjecture concerns
$V$-modules, not finitely generated $\N$-gradable weak $V$-modules.
Condition 3 would be a consequence of the finiteness of the dimensions
of genus-one conformal blocks, if the conformal field theory had been
constructed, and is thus natural and necessary. For vertex operator
algebras associated to affine Lie algebras (Wess-Zumino-Novikov-Witten
models) and vertex operator algebras associated to the Virasoro algebra
(minimal models), Condition 2 can be verified easily by reformulating
the corresponding complete reducibility results in terms of the
representation theory of affine Lie algebras and the Virasoro algebra.
For these vertex operator algebras, Condition 3 can also be easily
verified by using results in the representation theory of affine Lie
algebras and the Virasoro algebra.  In fact, Condition 3 was stated to
hold for these algebras in Zhu's paper \cite{Z} and was verified by
Dong-Li-Mason \cite{DLM} (see also \cite{AN} for the case of minimal
models). In addition, as we have mentioned in the abstract and 
introduction, by results of Li \cite{L} and 
Abe-Buhl-Dong \cite{ABD}, 
Conditions (ii) and (iii) are equivalent to a 
single condition that every weak $V$-module is completely 
reducible.

Because of this proposition, we see from Theorem 3.9
in \cite{H5} that $V$ satisfies all the conditions needed in the results
proved in \cite{HL1}--\cite{HL4}, \cite{H1} and \cite{HL6}. 
Thus we have:

\begin{thm}
Let $V$ be a vertex operator algebra satisfying the conditions above.
Then the category of $V$-modules has 
a natural structure of braided tensor category.
\end{thm}

The proofs of the rigidity and nondegeneracy in Sections 3 and 4,
respectively, depend not only on this theorem, but also
on the detailed construction of the braided tensor category
structure. For reader's convenience, we 
now briefly review the structures which are needed in the proof of the 
theorem above and in the main theorems of the present paper in these 
later sections.

Let $W_{1}$ and $W_{2}$ be $V$-modules.
In the present paper, we shall need the $P(z)$-tensor product
$W_{1}\boxtimes_{P(z)}W_{2}$ of $V$-modules $W_{1}$ and $W_{2}$ for
$z\in \mathbb{C}^{\times}$. 

The Jacobi identity for intertwining operators motivates a natural 
action $\tau_{P(z)}$ of 
$$x_{0}^{-1}\delta\left(\frac{x^{-1}_{1}-z}{x_{0}}
\right)
Y_{t}(v, x_{1})$$
on $(W_{1}\otimes W_{2})^{*}$ for $v\in V$ given by 
\begin{eqnarray*}
\lefteqn{\left(\tau_{P(z)}
\left(x_{0}^{-1}\delta\left(\frac{x^{-1}_{1}-z}{x_{0}}\right)
Y_{t}(v, x_{1})\right)\lambda\right)(w_{(1)}\otimes w_{(2)})=}\nonumber\\
&&=z^{-1}\delta\left(\frac{x^{-1}_{1}-x_{0}}{z}\right)
\lambda(Y_{1}(e^{x_{1}L(1)}(-x_{1}^{-2})^{L(0)}v, x_{0})w_{(1)}\otimes w_{(2)})
\nonumber\\
&&\quad +x^{-1}_{0}\delta\left(\frac{z-x^{-1}_{1}}{-x_{0}}\right)
\lambda(w_{(1)}\otimes Y_{2}^{*}(v, x_{1})w_{(2)})
\end{eqnarray*}
where $\lambda\in (W_{1}\otimes W_{2})^{*}$ and
$$Y_{t}(v, x_{1})=v\otimes x_{1}^{-1}\delta\left(\frac{t}{x_{1}}\right).$$
In particular, we have an action $\tau_{P(z)}(Y_{t}(v, x_{1}))$
of $Y_{t}(v, x_{1})$ on 
$(W_{1}\otimes W_{2})^{*}$.

Consider $\lambda\in (W_{1}\otimes W_{2})^{*}$ satisfying
the following conditions: 

\begin{description}
\item[The compatibility condition:]
{\bf (a)} The  {\it lower
truncation condition}:
For all $v\in V$, the formal Laurent series $\tau_{P(z)}(Y_{t}(v, x))\lambda$ 
involves only finitely many negative 
powers of $x$.

{\bf (b)} The following formula holds:
\begin{eqnarray*}
\lefteqn{\tau_{P(z)}\left(x_{0}^{-1}\delta\left(\frac{x^{-1}_{1}-z}{x_{0}}
\right)
Y_{t}(v, x_{1})\right)\lambda=}\nn
&&=x_{0}^{-1}\delta\left(\frac{x^{-1}_{1}-z}{x_{0}}\right)
\tau_{P(z)}(Y_{t}(v, x_{1}))\lambda
\end{eqnarray*}
for all $v\in V$.

\item[The local grading-restriction  condition:]
{\bf (a)} The {\it grading condition}: 
$\lambda$ is a (finite) sum of 
weight vectors (eigenvectors of $L(0)$) of $(W_{1}\otimes W_{2})^{*}$.

{\bf (b)} Let $W_{\lambda}$ be the smallest subspace of $(W_{1}\otimes
W_{2})^{*}$ containing $\lambda$ and stable under the component
operators $\tau_{P(z)}(v\otimes t^{n})$ of the operators 
$\tau_{P(z)}(Y_{t}(v,
x))$ for $v\in V$, $n\in \mathbb{Z}$. Then the weight spaces
$(W_{\lambda})_{(n)}$, $n\in \mathbb{C}$, of the (graded) space
$W_{\lambda}$ have the properties
\begin{eqnarray*}
&\mbox{\rm dim}\ (W_{\lambda})_{(n)}<\infty \;\;\;\mbox{\rm for}\
n\in \mathbb{C},&\\
&(W_{\lambda})_{(n)}=0 \;\;\;\mbox{\rm for $n$ such that $\Re{(n)}<<0$.}&
\end{eqnarray*}
\end{description}

Let $W_{1}\hboxtr_{P(z)} W_{2}$ be the subspace of 
$(W_{1}\otimes W_{2})^{*}$ consisting of all
elements satisfying the two conditions above.
Then by Theorem 13.5 in \cite{HL4} and Theorems 3.1 and 3.9 in  \cite{H5},
$W_{1}\hboxtr_{P(z)} W_{2}$ is a $V$-module.
The $P(z)$-tensor product module $W_{1}\boxtimes_{P(z)}W_{2}$ is 
defined to be the contragredient module 
$(W_{1}\hboxtr_{P(z)} W_{2})'$ of $W_{1}\hboxtr_{P(z)} W_{2}$.
We take the $P(1)$-tensor product $\boxtimes_{P(1)}$ to be the 
tensor product bifunctor and denote it by $\boxtimes$.

For any $V$-module
$W=\coprod_{n\in \mathbb{Q}}W_{(n)}$, we use 
$\overline{W}$ to denote its algebraic completion 
$\prod_{n\in \mathbb{Q}}W_{(n)}$. 
In addition to the $P(z)$-tensor product module, by the definition of 
$P(z)$-tensor product (Definition 4.1 in \cite{HL1})
and Proposition 12.2
in \cite{HL4},
there is also 
an intertwining operator $\Y$ of type ${W_{1}\boxtimes_{P(z)}W_{2}\choose
W_{1}W_{2}}$ such that the following universal property holds:
For any intertwining operator $\tilde{\Y}$ of type 
${W_{3}\choose
W_{1}W_{2}}$, there is a module map $f: W_{1}\boxtimes_{P(z)}W_{2}
\to W_{3}$ such that $\tilde{\Y}=f\circ \Y$. For $w_{1}\in W_{1}$ and 
$w_{2}\in W_{2}$, the $P(z)$-tensor product of elements
$w_{1}$ and $w_{2}$
to be 
$$w_{1}\boxtimes w_{2}=\Y(w_{1}, z)w_{2}\in 
\overline{W_{1}\boxtimes_{P(z)}W_{2}},$$
where we use the convention 
$$\Y(w_{1}, z)w_{2}=\Y(w_{1}, x)w_{2}|_{x^{n}=e^{n\log z},\;n\in \C}$$
and 
$$\log z=\log |z|+i\arg z,\; 0\le \arg z<2\pi.$$
(We shall continue to use this convention through this paper.)

The existence of these tensor products of elements is a very important 
feature of the tensor product theory. They provide a
powerful tool for proving theorems in the tensor product theory:
We first prove the results on these tensor products of elements. Since 
the homogeneous components of these tensor products of elements span
the tensor product modules, we obtain the results we are interested. 
(One subtle thing is worth mentioning here. The space spanned by all tensor 
products of elements has almost no intersection with the tensor product 
module; in general, only elements of the form 
$\mathbf{1}\boxtimes_{P(z)} w$ are in $V\boxtimes_{P(z)}W$ where $W$ is 
an arbitrary $V$-module.) 

For $P(z)$-tensor products with different $z$, we have 
parallel transport isomorphisms between them. In this paper, 
as we have used above, for $z\in 
\C^{\times}$,
we shall use $\log z$ to denote that value of the logarithm of 
$z$ satisfying $0\le \Im(\log z)<2\pi$.
Let $W_{1}$ and $W_{2}$ be $V$-modules and $z_{1}, z_{2}\in \C^{\times}$. 
Giving a path $\gamma$ in
$\mathbb{C}^{\times}$ from $z_{1}$ to $z_{2}$.
The parallel isomorphism 
$\mathcal{T}_{\gamma}: W_{1}\boxtimes_{P(z_{1})}W_{2}\to 
W_{1}\boxtimes_{P(z_{2})}W_{2}$ (see \cite{HL3} \cite{HL6}) 
is given as follows: 
Let $\mathcal{Y}$ be the intertwining 
operator associated to  the $P(z_{2})$-tensor product 
$W_{1}\boxtimes_{P(z_{2})}W_{2}$
and $l(z_{1})$ the value of the logarithm of $z_{1}$ determined uniquely by
$\log z_{2}$ and the 
path $\gamma$. Then  $\mathcal{T}_{\gamma}$ is characterized by
$$\overline{\mathcal{T}}_{\gamma}(w_{1}
\boxtimes_{P(z_{1})}w_{2})=\mathcal{Y}(w_{1}, x)w_{2}
|_{x^{n}=e^{nl(z_{1})}, \; n\in \C}$$
for $w_{1}\in W_{1}$ and $w_{2}\in W_{2}$, where 
$\overline{\mathcal{T}}_{\gamma}$
is the natural extension of $\mathcal{T}_{\gamma}$ 
to the algebraic 
completion $\overline{W_{1}\boxtimes_{P(z_{1})} W_{2}}$ 
of $W_{1}\boxtimes_{P(z_{1})} W_{2}$.  The parallel isomorphism 
depends only on the homotopy class of $\gamma$.

For $z\in \C^{\times}$, 
the commutativity isomorphism for the $P(z)$-tensor product
is characterized as follows: Let $\gamma_{z}^{-}$ be a path
from $-z$ to $z$  in the 
closed upper half plane with $0$ deleted and
$\mathcal{T}_{\gamma_{z}^{-}}$ the corresponding 
parallel transport isomorphism. Then
$$\overline{\mathcal{C}_{P(z)}}(w_{1}\boxtimes_{P(z)} w_{2})=e^{zL(-1)}
\overline{\mathcal{T}}_{\gamma_{z}^{-}}
(w_{2}\boxtimes_{P(z)} w_{1})$$
where $w_{1}\in W_{1}$, $w_{2}\in W_{2}$. When $z=1$, we obtain 
a commutativity isomorphism 
$$\mathcal{C}_{P(1)}: W_{1}\boxtimes W_{2}\to 
W_{2}\boxtimes W_{1}$$
and we shall denote it simply as $\mathcal{C}$.

The tensor product  $w_{1}\boxtimes w_{2}$ of 
$w_{1}\in W_{1}$ and $w_{2}\in W_{2}$ is obtained from 
two elements $w_{1}$ and $w_{2}$ of the $V$-modules $W_{1}$ and $W_{2}$,
respectively. We also need
tensor products of more than two elements. Here we describe 
tensor products of three elements briefly. 
Let $W_{1}, W_{2}, W_{3}$ be $V$-modules and $w_{1}\in W_{1}$,
$w_{2}\in W_{2}$ and $w_{3}\in W_{3}$. Let $z_{1}$ and 
$z_{2}$ be two nonzero complex numbers.
Since in general 
$w_{2}\boxtimes_{P(z_{2})} w_{3}$ does not belong to 
$W_{2}\boxtimes_{P(z_{2})}W_{3}$, we cannot define 
$w_{1}\boxtimes_{P(z_{1})}(w_{2}\boxtimes_{P(z_{2})} w_{3})$
simply to be the $P(z_{2})$-tensor product of $w_{1}$ and 
$w_{2}\boxtimes_{P(z_{2})} w_{3}$.  But using the convergence
of products of intertwining operators, the series 
$$\sum_{n\in \mathbb{Z}}
w_{1}\boxtimes_{P(z_{1})}P_{n}(w_{2}\boxtimes_{P(z_{2})} w_{3})$$
($P_{n}$
is the projection map from a $V$-module to the subspace of weight $n$)
is absolutely 
convergent in a natural sense when $|z_{1}|>|z_{2}|$ and the sum is 
in $\overline{W_{1}\boxtimes_{P(z_{1})}(W_{2}
\boxtimes_{P(z_{2})} W_{3})}$.
We define
$w_{1}\boxtimes_{P(z_{1})}(w_{2}\boxtimes_{P(z_{2})} w_{3})$
to be this sum.
Similarly, when $|z_{2}|>|z_{1}-z_{2}|>0$, we have 
$$(w_{1}\boxtimes_{P(z_{1}-z_{2})}w_{2})\boxtimes_{P(z_{2})} w_{3}\in
\overline{(W_{1}\boxtimes_{P(z_{1}-z_{2})}W_{2})\boxtimes_{P(z_{2})} 
W_{3}}.$$
The homogeneous components of these 
tensor products of three elements also
span the corresponding tensor products of three modules.

Let $z_{1}, z_{2}$ be complex numbers satisfying 
$|z_{1}|>|z_{2}|>|z_{1}-z_{2}|>0$ and
$W_{1}$, $W_{2}$ and $W_{3}$ $V$-modules. Then 
an associativity isomorphism 
$$\A_{P(z_{1}), P(z_{2})}^{P(z_{1}-z_{2}), P(z_{2})}:
W_{1}\boxtimes_{P(z_{1})}(W_{2}\boxtimes_{P(z_{2})}W_{3})
\to (W_{1}\boxtimes_{P(z_{1}-z_{2})}W_{2})\boxtimes_{P(z_{2})}W_{3}$$
was constructed by the author in \cite{H1} and \cite{H5} and  is 
characterized by the property
$$\overline{\A_{P(z_{1}, P(z_{2})}^{P(z_{1}-z_{2}), P(z_{2})}}
(w_{1}\boxtimes_{P(z_{1})}(w_{2}\boxtimes_{P(z_{2})}w_{3}))
=(w_{1}\boxtimes_{P(z_{1}-z_{2})}w_{2})\boxtimes_{P(z_{2})}w_{3}$$
for $w_{1}\in W_{1}$, $w_{2}\in W_{2}$ and $w_{3}\in W_{3}$, 
where 
$$\overline{\A_{P(z_{1}, P(z_{2})}^{P(z_{1}-z_{2}), P(z_{2})}}: 
\overline{W_{1}\boxtimes_{P(z_{1})}(W_{2}\boxtimes_{P(z_{2})}W_{3})}
\to \overline{(W_{1}\boxtimes_{P(z_{1}-z_{2})}W_{2})
\boxtimes_{P(z_{2})}W_{3}}$$
is the natural extension of 
$\A_{P(z_{1}), P(z_{2})}^{P(z_{1}-z_{2}), P(z_{2})}$ to the 
algebraic completion 
of $W_{1}\boxtimes_{P(z_{1})}(W_{2}\boxtimes_{P(z_{2})}W_{3})$.

To obtain the associativity isomorphism 
$$\mathcal{A}: 
W_{1}\boxtimes (W_{2}\boxtimes W_{3})\to 
(W_{1}\boxtimes W_{2})\boxtimes W_{3}$$
for the braided tensor category structure,
we need certain parallel isomorphisms. Let 
$z_{1}$ and $z_{2}$ be real numbers satisfying 
$z_{1}>z_{2}>z_{1}-z_{2}\ge 0$.  Let $\gamma_{1}$ and $\gamma_{2}$ 
be paths in $(0, \infty)$
from $1$ to $z_{1}$ and $z_{2}$, respectively, and 
$\gamma_{3}$ and $\gamma_{4}$ be 
paths in $(0, \infty)$
from  $z_{2}$ and $z_{1}-z_{2}$ to $1$, respectively. 
Then the associativity isomorphism for the braided tensor category 
structure on the module category for $V$ is given by 
$$\mathcal{A}=\mathcal{T}_{\gamma_{3}}\circ (\mathcal{T}_{\gamma_{4}}
\boxtimes_{P(z_{2})} I_{W_{3}})\circ 
\mathcal{A}^{P(z_{1}-z_{2}), P(z_{2})}_{P(z_{1}), P(z_{2})}\circ
(I_{W_{1}} \boxtimes_{P(z_{1})} 
\mathcal{T}_{\gamma_{2}})\circ \mathcal{T}_{\gamma_{1}},$$
that is, given by the commutative diagram
$$
\begin{CD}
W_{1}\boxtimes_{P(z_{1})} (W_{2}
\boxtimes_{P(z_{2})} W_{3})
@>\mathcal{A}_{P(z_{2}), P(z_{2})}^{P(z_{2}-z_{3}), P(z_{3})}>>
(W_{1}\boxtimes_{P(z_{1}-z_{2})} W_{2})
\boxtimes_{P(z_{2})} W_{3}\\ 
@A(I_{W_{1}} \boxtimes_{P(z_{1})} 
\mathcal{T}_{\gamma_{2}})\circ \mathcal{T}_{\gamma_{1}}AA 
@VV\mathcal{T}_{\gamma_{3}}\circ (\mathcal{T}_{\gamma_{4}}
\boxtimes_{P(z_{2})} I_{W_{3}})V\\
W_{1}\boxtimes (W_{2}
\boxtimes W_{3}) @>\mathcal{A}>> 
(W_{1}\boxtimes W_{2})
\boxtimes W_{3}
\end{CD}
$$

The coherence properties are now easy consequences of
the constructions and characterizations of the associativity 
and commutativity isomorphisms. Here we sketch the 
proof of the commutativity of the pentagon diagram. 
Let $W_{1}$, $W_{2}$, $W_{3}$ and $W_{4}$ be $V$-modules and let $z_{1},
z_{2}, z_{3}\in \C$ satisfying 
\begin{eqnarray*}
&|z_{1}|>|z_{2}|>|z_{3}|>
|z_{1}-z_{3}|>|z_{2}-z_{3}|>|z_{1}-z_{2}|>0,&\\
&|z_{1}|>|z_{2}-z_{3}|+|z_{3}|,&\\
&|z_{2}|>|z_{1}-z_{2}|+|z_{3}|,&\\
&|z_{2}|>|z_{1}-z_{2}|+|z_{2}-z_{3}|.&
\end{eqnarray*}
For example, we can take $z_{1}=7$, $z_{2}=6$ and $z_{3}=4$. 
We first prove the commutativity of the following diagram:

\vspace{3.5em}

\begin{picture}(150,100)(-70,0)

\put(-85,20){\footnotesize $((W_{1}\boxtimes_{P(z_{12})} 
W_{2})\boxtimes_{P(z_{23})} W_{3})
\boxtimes_{P(z_{3})} W_{4}$}       
\put(123,20){\footnotesize $(W_{1}\boxtimes_{P(z_{13})} 
(W_{2}\boxtimes_{P(z_{23})} W_{3}))
\boxtimes_{P(z_{3})} W_{4}$}         
\put(-85,68){\footnotesize $(W_{1}\boxtimes_{P(z_{12})} 
W_{2})\boxtimes_{P(z_{2})} (W_{3}
\boxtimes_{P(z_{3})} W_{4})$}    
\put(123,68){\footnotesize $W_{1}\boxtimes_{P(z_{1})} 
((W_{2}\boxtimes_{P(z_{23})} W_{3})
\boxtimes_{P(z_{3})} W_{4}))$} 
\put(20,116){\footnotesize $W_{1}\boxtimes_{P(z_{1})} 
(W_{2}\boxtimes_{P(z_{2})} (W_{3}
\boxtimes_{P(z_{3})} W_{4}))$}          

\put(120,23){\vector(-1,0){25}}  
\put(5,60){\vector(0,-1){28}} 
\put(212,60){\vector(0,-1){28}} 
 
\put(100,105){\vector(-3,-1){75}} 
\put(110,105){\vector(3,-1){75}}
\end{picture}
\begin{equation}\label{pent1}
\end{equation}
where $z_{12}=z_{1}-z_{2}$ and $z_{23}=z_{2}-z_{3}$.
For $w_{1}\in W_{1}$, $w_{2}\in W_{2}$, $w_{3}\in W_{3}$ and 
$w_{4}\in W_{4}$, we consider 
$$w_{1}\boxtimes_{P(z_{1})} (w_{2}\boxtimes_{P(z_{2})}
(w_{3}\boxtimes_{P(z_{3})}w_{4}))\in 
\overline{W_{1}\boxtimes_{P(z_{1})} (W_{2}\boxtimes_{P(z_{2})}
(W_{3}\boxtimes_{P(z_{3})}W_{4}))}.$$
By the characterizations of the associativity isomorphisms, 
we see that the compositions of the natural extensions 
of the module maps in the two routes in (\ref{pent1}) applying to 
this element both give 
$$((w_{1}\boxtimes_{P(z_{12})} w_{2})\boxtimes_{P(z_{23})}
w_{3})\boxtimes_{P(z_{3})}w_{4}\in
\overline{((W_{1}\boxtimes_{P(z_{12})} W_{2})\boxtimes_{P(z_{23})}
W_{3})\boxtimes_{P(z_{3})}W_{4}}.$$
Since the homogeneous components of 
$$w_{1}\boxtimes_{P(z_{1})} (w_{2}\boxtimes_{P(z_{2})}
(w_{3}\boxtimes_{P(z_{3})}w_{4}))$$
for $w_{1}\in W_{1}$, $w_{2}\in W_{2}$, $w_{3}\in W_{3}$ and 
$w_{4}\in W_{4}$ span 
$$W_{1}\boxtimes_{P(z_{1})} (W_{2}\boxtimes_{P(z_{2})}
(W_{3}\boxtimes_{P(z_{3})}W_{4})),$$
the diagram (\ref{pent1}) above is commutative. 

On the other hand, by the definition of $\A$, the  diagrams 
\begin{equation}\label{pent2}
\begin{picture}(60,90)(20,0)
\put(-145,68){\footnotesize $W_{1}\boxtimes_{P(z_{1})} 
(W_{2}\boxtimes_{P(z_{2})} (W_{3}
\boxtimes_{P(z_{3})} W_{4}))$}
\put(63,68){\footnotesize $(W_{1}\boxtimes_{P(z_{12})} 
W_{2})\boxtimes_{P(z_{2})} (W_{3}
\boxtimes_{P(z_{3})} W_{4})$} 
\put(-105,20){\footnotesize $W_{1}\boxtimes
(W_{2}\boxtimes (W_{3}
\boxtimes W_{4}))$}
\put(88,20){\footnotesize $(W_{1}\boxtimes
W_{2})\boxtimes (W_{3}
\boxtimes W_{4})$}

\put(5,23){\vector(1,0){78}}  
\put(30,71){\vector(1,0){25}}  
\put(-55,60){\vector(0,-1){28}} 
\put(152,60){\vector(0,-1){28}} 
\end{picture}
\end{equation}
\begin{equation}\label{pent3}
\begin{picture}(60,90)(20,0)
\put(-145,68){\footnotesize $(W_{1}\boxtimes_{P(z_{12})} 
W_{2})\boxtimes_{P(z_{2})} (W_{3}
\boxtimes_{P(z_{3})} W_{4})$}
\put(63,68){\footnotesize $((W_{1}\boxtimes_{P(z_{12})} 
W_{2})\boxtimes_{P(z_{23})} W_{3})
\boxtimes_{P(z_{3})} W_{4}$}
\put(-105,20){\footnotesize $(W_{1}\boxtimes
W_{2})\boxtimes (W_{3}
\boxtimes W_{4})$}
\put(88,20){\footnotesize $((W_{1}\boxtimes
W_{2})\boxtimes W_{3})
\boxtimes W_{4}$}

\put(5,23){\vector(1,0){78}}  
\put(34,71){\vector(1,0){25}}
\put(-55,60){\vector(0,-1){28}} 
\put(152,60){\vector(0,-1){28}} 
\end{picture}
\end{equation}
\begin{equation}\label{pent4}
\begin{picture}(60,90)(20,0)
\put(-145,68){\footnotesize $W_{1}\boxtimes_{P(z_{1})} 
(W_{2}\boxtimes_{P(z_{2})} (W_{3}
\boxtimes_{P(z_{3})} W_{4}))$}
\put(63,68){\footnotesize $W_{1}\boxtimes_{P(z_{1})} 
((W_{2}\boxtimes_{P(z_{23})} W_{3})
\boxtimes_{P(z_{3})} W_{4}))$} 
\put(-105,20){\footnotesize $W_{1}\boxtimes
(W_{2}\boxtimes (W_{3}
\boxtimes W_{4}$))} 
\put(88,20){\footnotesize $W_{1}\boxtimes
((W_{2}\boxtimes W_{3})
\boxtimes W_{4})$} 

\put(5,23){\vector(1,0){78}}  
\put(34,71){\vector(1,0){25}}
\put(-55,60){\vector(0,-1){28}} 
\put(152,60){\vector(0,-1){28}} 
\end{picture}
\end{equation}
\begin{equation}\label{pent5}
\begin{picture}(60,80)(20,0)
\put(-145,68){\footnotesize $W_{1}\boxtimes_{P(z_{1})} 
((W_{2}\boxtimes_{P(z_{23})} W_{3})
\boxtimes_{P(z_{3})} W_{4}))$} 
\put(63,68){\footnotesize $(W_{1}\boxtimes_{P(z_{13})} 
(W_{2}\boxtimes_{P(z_{23})} W_{3}))
\boxtimes_{P(z_{3})} W_{4}$}
\put(-105,20){\footnotesize $W_{1}\boxtimes
((W_{2}\boxtimes W_{3})
\boxtimes W_{4})$}
\put(88,20){\footnotesize $(W_{1}\boxtimes
(W_{2}\boxtimes W_{3}))
\boxtimes W_{4}$} 

\put(5,23){\vector(1,0){78}}  
\put(34,71){\vector(1,0){25}}
\put(-55,60){\vector(0,-1){28}} 
\put(152,60){\vector(0,-1){28}} 
\end{picture}
\end{equation}
\begin{equation}\label{pent6}
\begin{picture}(60,90)(20,0)
\put(-145,68){\footnotesize $(W_{1}\boxtimes_{P(z_{13})} 
(W_{2}\boxtimes_{P(z_{23})} W_{3}))
\boxtimes_{P(z_{3})} W_{4}$}
\put(63,68){\footnotesize $((W_{1}\boxtimes_{P(z_{12})} 
W_{2})\boxtimes_{P(z_{23})} W_{3})
\boxtimes_{P(z_{3})} W_{4}$}   
\put(-105,20){\footnotesize $(W_{1}\boxtimes
(W_{2}\boxtimes W_{3}))
\boxtimes W_{4}$} 
\put(88,20){\footnotesize $((W_{1}\boxtimes
W_{2})\boxtimes W_{3})
\boxtimes W_{4}$}      

\put(5,23){\vector(1,0){78}}  
\put(34,71){\vector(1,0){25}}
\put(-55,60){\vector(0,-1){28}} 
\put(152,60){\vector(0,-1){28}} 
\end{picture}
\end{equation}
are all commutative. Combining all the diagrams 
(\ref{pent1})--(\ref{pent6}) above, we see that
the pentagon diagram

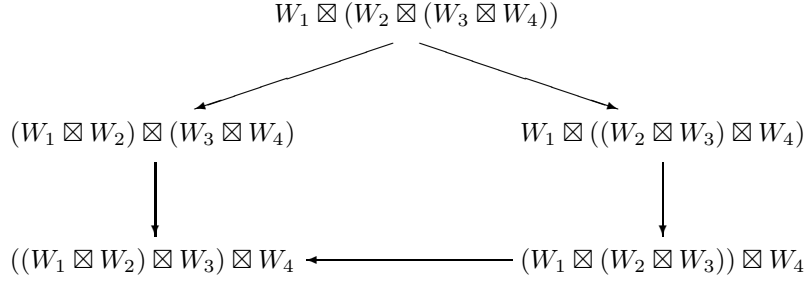
\vspace{3.5em}

\begin{picture}(150,100)(-70,0)

\put(-45,20){\footnotesize $((W_{1}\boxtimes
W_{2})\boxtimes W_{3})
\boxtimes W_{4}$}       
\put(148,20){\footnotesize $(W_{1}\boxtimes
(W_{2}\boxtimes W_{3}))
\boxtimes W_{4}$}         
\put(-45,68){\footnotesize $(W_{1}\boxtimes
W_{2})\boxtimes (W_{3}
\boxtimes W_{4})$}    
\put(148,68){\footnotesize $W_{1}\boxtimes
((W_{2}\boxtimes W_{3})
\boxtimes W_{4})$} 
\put(55,113){\footnotesize $W_{1}\boxtimes
(W_{2}\boxtimes (W_{3}
\boxtimes W_{4}))$}          

\put(145,23){\vector(-1,0){78}}  
\put(10,60){\vector(0,-1){28}} 
\put(202,60){\vector(0,-1){28}} 
 
\put(100,105){\vector(-3,-1){75}} 
\put(110,105){\vector(3,-1){75}}
\end{picture}

\noindent is also commutative. 

The proof of the commutativity of the hexagon diagrams
is similar.

The unit object is $V$. For any $z\in \C^{\times}$ and any $V$-module $W$, 
the left $P(z)$-unit isomorphism $l_{W; z}: V\boxtimes_{P(z)}W
\to W$ is characterized by 
$$l_{W; z}(\mathbf{1}\boxtimes_{P(z)} w)=w$$
for $w\in W$ and the right 
$P(z)$-unit isomorphism $r_{W; z}: W\boxtimes_{P(z)}V
\to W$ is characterized by 
$$\overline{r_{W; z}}(w\boxtimes_{P(z)} \mathbf{1})=e^{zL(-1)} w$$
for $w\in W$. 
In particular, we have the 
left unit isomorphism $l_{W}=l_{W; 1}: V\boxtimes W
\to W$ and the right unit isomorphism $r_{W}=r_{W; 1}: W\boxtimes V
\to W$. The proof of the commutativity of the 
diagrams for unit isomorphisms is similar to the above 
proof of the commutativity of the pentagon diagram. 

\renewcommand{\theequation}{\thesection.\arabic{equation}}
\renewcommand{\thethm}{\thesection.\arabic{thm}}
\setcounter{equation}{0}
\setcounter{thm}{0}

\section{The fusing and braiding matrices, modular 
transformations and the Verlinde 
conjecture}

In the proofs in the next two sections 
of the rigidity and nondegeneracy of the 
semisimple braided tensor category of $V$-modules for a vertex operator 
algebra $V$ satisfying the conditions in the preceding section, 
we need fusing and braiding matrices and 
certain consequences of the Verlinde conjecture, which
was recently proved by the author \cite{H7} for a vertex operator 
algebra satisfying conditions slightly weaker than those assumed in
the present paper. Here we give a brief review. For details, see 
\cite{H4}--\cite{H7}.

Using the theory of associative algebras and 
Zhu's  algebra for a vertex operator algebra,
it is easy to see  that for a vertex operator algebra $V$ 
such that every $\N$-gradable $V$-module is completely reducible,
there are only finitely many
inequivalent irreducible $V$-modules (see Theorem 3.2 in \cite{DLM}).  
Let $\A$ be the set of
equivalence classes of irreducible $V$-modules. We denote the
equivalence class containing $V$ by $e$. By the main result of 
\cite{AM} and Theorem 11.3 in \cite{DLM},  we know that 
$V$-modules are all graded by $\R$. 
For each $a\in \A$,
we choose a representative $W^{a}$ of $a$ such that 
$W^{e}=V$. For $a\in \A$, let 
$h_{a}$ be the lowest weight of $W^{a}$, that is, $h_{a}\in \R$
such that $W^{a}=\coprod_{n\in h_{a}+\N}W^{a}_{(n)}$.
By Propositions 5.3.1 and 5.3.2 in \cite{FHL}, the contragredient
module of an irreducible module is also irreducible and the
contragredient module of the contragredient module of a $V$-module 
is naturally equivalent to the $V$-module itself.
So we have a bijective map
\begin{eqnarray*}
^{\prime}:
\mathcal{A}&\to& \mathcal{A}\\
a&\mapsto& a'.
\end{eqnarray*}

Let $\mathcal{V}_{a_{1}a_{2}}^{a_{3}}$ for $a_{1}, a_{2}, a_{3}
\in \mathcal{A}$
be the space of intertwining operators of 
type ${W^{a_{3}}\choose W^{a_{1}}W^{a_{2}}}$ and 
$N_{a_{1}a_{2}}^{a_{3}}$ for $a_{1}, a_{2}, a_{3}
\in \mathcal{A}$ the fusion rule, that is, the dimension of 
the space of intertwining operators of 
type ${W^{a_{3}}\choose W^{a_{1}}W^{a_{2}}}$. 
The fusion rules 
$N_{a_{1}a_{2}}^{a_{3}}$ for 
$a_{1}, a_{2}, a_{3}\in \mathcal{A}$ are all finite \cite{GN}
\cite{L} \cite{AN} \cite{H5}.

We now discuss matrix elements of fusing and braiding isomorphisms. 
We need 
to use different bases of one space of intertwining 
operators. We shall  use $p=1, 2, 3, 4, 5, 6, \dots$
to label different bases. 
For  $a_{1}, a_{2}, 
a_{3}\in \mathcal{A}$ and $p=1, 2, 3, 4, 5, 6, \dots$, let 
$\{\Y_{a_{1}a_{2}; i}^{a_{3}; (p)}\;|\; i=1, \dots, 
N_{a_{1}a_{2}}^{a_{3}}\}$, be 
bases of $\mathcal{V}_{a_{1}a_{2}}^{a_{3}}$. 
The associativity of intertwining operators proved and studied in 
\cite{H1}, \cite{H4}
and \cite{H5} says that 
there exist 
$$F(\Y_{a_{1}a_{5}; i}^{a_{4}; (1)}\otimes 
\Y_{a_{2}a_{3}; j}^{a_{5}; (2)}; 
\Y_{a_{6}a_{3}; l}^{a_{4}; (3)}
\otimes \Y_{a_{1}a_{2}; k}^{a_{6}; (4)}) \in \C$$
for $a_{1}, \dots, a_{6}\in \mathcal{A}$, $i=1, \dots, 
N_{a_{1}a_{5}}^{a_{4}}$, $j=1, \dots, 
N_{a_{2}a_{3}}^{a_{5}}$, $k=1, \dots, 
N_{a_{6}a_{3}}^{a_{4}}$, $l=1, \dots, 
N_{a_{1}a_{2}}^{a_{6}}$
such that 
\begin{eqnarray*}
\lefteqn{\langle w'_{a_{4}}, \Y_{a_{1}a_{5}; i}^{a_{4}; (1)}(w_{a_{1}}, z_{1})
\Y_{a_{2}a_{3}; j}^{a_{5}; (2)}(w_{a_{2}}, z_{2})w_{a_{3}}\rangle}\nn
&&=\sum_{a_{6}\in \A}
\sum_{k=1}^{N_{a_{6}a_{3}}^{a_{4}}}\sum_{l=1}^{N_{a_{1}a_{2}}^{a_{6}}}
F(\Y_{a_{1}a_{5}; i}^{a_{4}; (1)}\otimes \Y_{a_{2}a_{3}; j}^{a_{5}; (2)}; 
\Y_{a_{6}a_{3}; l}^{a_{4}; (3)}
\otimes \Y_{a_{1}a_{2}; k}^{a_{6}; (4)})\cdot\nn
&&\quad \quad \quad \quad \quad \quad \quad \quad \quad \cdot 
\langle w'_{a_{4}}, 
\Y_{a_{6}a_{3}; k}^{a_{4}; (3)}(\Y_{a_{1}a_{2}; l}^{a_{6}; (4)}(w_{a_{1}}, z_{1}-z_{2})
w_{a_{2}}, z_{2})w_{a_{3}}\rangle
\nn
&&
\end{eqnarray*}
when $|z_{1}|>|z_{2}|>|z_{1}-z_{2}|>0$, 
for $a_{1}, \dots, a_{5}\in \A$, $w_{a_{1}}\in W^{a_{1}}$, 
$w_{a_{2}}\in W^{a_{2}}$, $w_{a_{3}}\in W^{a_{3}}$,  
$w'_{a_{4}}\in (W^{a_{4}})'$, $i=1, \dots, 
N_{a_{1}a_{5}}^{a_{4}}$ and $j=1, \dots, 
N_{a_{2}a_{3}}^{a_{5}}$. The numbers 
$$F(\Y_{a_{1}a_{5}; i}^{a_{4}; (1)}\otimes \Y_{a_{2}a_{3}; j}^{a_{5}; (2)}; 
\Y_{a_{6}a_{3}; k}^{a_{4}; (3)}
\otimes \Y_{a_{1}a_{2}; l}^{a_{6}; (4)})$$
together give a matrix which represents a linear 
isomorphism
$$\coprod_{a_{1}, a_{2}, a_{3}, a_{4}, a_{5}\in \A}
\mathcal{V}_{a_{1}a_{5}}^{a_{4}}\otimes 
\mathcal{V}_{a_{2}a_{3}}^{a_{5}}\to 
\coprod_{a_{1}, a_{2}, a_{3}, a_{4}, a_{6}\in \A}
\mathcal{V}_{a_{6}a_{3}}^{a_{4}}
\otimes \mathcal{V}_{a_{1}a_{2}}^{a_{6}},$$
called the {\it fusing isomorphism}, 
such that these numbers 
are the matrix elements.

By the commutativity of intertwining operators proved 
and studied in \cite{H2}, \cite{H4} and \cite{H5},
$r\in \Z$, 
there exist 
$$B^{(r)}(\Y_{a_{1}a_{5}; i}^{a_{4}; (1)}
\otimes \Y_{a_{2}a_{3}; j}^{a_{5}; (2)}; 
\Y_{a_{2}a_{6}; l}^{a_{4}; (3)}
\otimes \Y_{a_{1}a_{3}; k}^{a_{6}; (4)}) \in \C$$
for $a_{1}, \dots, a_{6}\in \mathcal{A}$, $i=1, \dots, 
N_{a_{1}a_{5}}^{a_{4}}$, $j=1, \dots, 
N_{a_{2}a_{3}}^{a_{5}}$, $k=1, \dots, 
N_{a_{2}a_{6}}^{a_{4}}$, $l=1, \dots, 
N_{a_{1}a_{3}}^{a_{6}}$, 
such that 
the analytic extension of the single-valued analytic function
$$\langle w'_{a_{4}}, \Y_{a_{1}a_{5}; i}^{a_{4}; (1)}(w_{a_{1}}, z_{1})
\Y_{a_{2}a_{3}; j}^{a_{5}; (2)}(w_{a_{2}}, z_{2})w_{a_{3}}\rangle$$
on the region $|z_{1}|>|z_{2}|>0$, $0\le \arg z_{1}, \arg z_{2}<2\pi$
along the path 
$$t \mapsto \left(\frac{3}{2}
-\frac{e^{(2r+1)\pi i t}}{2}, \frac{3}{2}
+\frac{e^{(2r+1)\pi i t}}{2}\right)$$ 
to the region $|z_{2}|>|z_{1}|>0$, $0\le \arg z_{1}, \arg z_{2}<2\pi$
is 
\begin{eqnarray*}
\lefteqn{\sum_{a_{6}\in \A}
\sum_{k=1}^{N_{a_{2}a_{6}}^{a_{4}}}\sum_{l=1}^{N_{a_{1}a_{3}}^{a_{6}}}
B^{(r)}(\Y_{a_{1}a_{5}; i}^{a_{4}; (1)}\otimes \Y_{a_{2}a_{3}; j}^{a_{5}; (2)};
\Y_{a_{2}a_{6}; k}^{a_{4}; (3)}\otimes 
\Y_{a_{1}a_{3}; l}^{a_{6}; (4)})\cdot}\nn
&&\quad \quad \quad \quad \quad \quad \quad \quad \quad \cdot 
\langle w'_{a_{4}}, 
\Y_{a_{2}a_{6}; k}^{a_{4}; (3)}(
w_{a_{2}}, z_{1})\Y_{a_{1}a_{3}; l}^{a_{6}; (4)}(w_{a_{1}}, z_{2})
w_{a_{3}}\rangle.
\end{eqnarray*}
The numbers 
$$B^{(r)}(\Y_{a_{1}a_{5}; i}^{a_{4}; (1)}\otimes \Y_{a_{2}a_{3}; j}^{a_{5}; (2)};
\Y_{a_{2}a_{6}; k}^{a_{4}; (3)}\otimes \Y_{a_{1}a_{3}; l}^{a_{6}; (4)})$$
together give a linear isomorphism 
$$B^{(r)}: \coprod_{a_{1}, a_{2}, a_{3}, a_{4}, a_{5}\in \A}
\mathcal{V}_{a_{1}a_{5}}^{a_{4}}\otimes 
\mathcal{V}_{a_{2}a_{3}}^{a_{5}}\to 
\coprod_{a_{1}, a_{2}, a_{3}, a_{4}, a_{6}\in \A}
\mathcal{V}_{a_{2}a_{6}}^{a_{4}}
\otimes \mathcal{V}_{a_{1}a_{3}}^{a_{6}},$$
called the {\it braiding isomorphism}, 
such that these numbers 
are the matrix elements.

In this paper, we are mainly interested in the square 
$(B^{(r)})^{2}$ of $B^{(r)}$. We shall also use similar 
notations to denote the matrix elements of
the square $(B^{(r)})^{2}$ of $B^{(r)}$ under the bases above as 
$$(B^{(r)})^{2}(\Y_{a_{1}a_{5}; i}^{a_{4}; (1)}\otimes \Y_{a_{2}a_{3}; j}^{a_{5}; (2)};
\Y_{a_{1}a_{6}; i}^{a_{4}; (3)}\otimes \Y_{a_{2}a_{3}; j}^{a_{6}; (4)}).$$
Let
\begin{eqnarray*}
\lefteqn{(B^{(r)})^{2}(\langle w_{a'_{4}}, \Y_{a_{1}a_{5}; i}^{a_{4}; (1)}(w_{a_{1}}, z_{1})
\Y_{a_{2}a_{3}; j}^{a_{5}; (2)}(w_{a_{2}}, z_{2})w_{a_{3}}\rangle)}\nn
&&=\sum_{a_{6}\in \A}
\sum_{k=1}^{N_{a_{2}a_{6}}^{a_{4}}}\sum_{l=1}^{N_{a_{1}a_{3}}^{a_{6}}}
(B^{(r)})^{2}(\Y_{a_{1}a_{5}; i}^{a_{4}; (1)}\otimes \Y_{a_{2}a_{3}; j}^{a_{5}; (2)};
\Y_{a_{1}a_{6}; i}^{a_{4}; (3)}\otimes \Y_{a_{2}a_{3}; j}^{a_{6}; (4)})\cdot\nn
&&\quad\quad\quad\quad\cdot 
\langle w_{a'_{4}}, \Y_{a_{1}a_{6}; i}^{a_{4}; (1)}(w_{a_{1}}, z_{1})
\Y_{a_{2}a_{3}; j}^{a_{6}; (2)}(w_{a_{2}}, z_{2})w_{a_{3}}\rangle.
\end{eqnarray*}
Then by definition, it is in fact the 
monodromy of the multi-valued analytic extension of 
$$\langle w_{a'_{4}}, \Y_{a_{1}a_{5}; i}^{a_{4}; (1)}(w_{a_{1}}, z_{1})
\Y_{a_{2}a_{3}; j}^{a_{5}; (2)}(w_{a_{2}}, z_{2})w_{a_{3}}\rangle$$
from $(z_{1}, z_{2})$ in the region $|z_{1}|>|z_{2}|>0$ 
to itself along the product of the path in the definition of 
$B^{(r)}$ above with itself (see Section 1 in \cite{H7} for more details).

We need an action of $S_{3}$ on the space 
$$\mathcal{V}=\coprod_{a_{1}, a_{2}, a_{3}\in 
\mathcal{A}}\mathcal{V}_{a_{1}a_{2}}^{a_{3}}.$$
For  $r\in \Z$, $a_{1}, a_{2}, a_{3}\in \mathcal{A}$, we have 
skew-symmetry isomorphisms
$\Omega_{-r}: \mathcal{V}_{a_{1}a_{2}}^{a_{3}} \to 
\mathcal{V}_{a_{2}a_{1}}^{a_{3}}$ and contragredient isomorphisms
$A_{-r}: \mathcal{V}_{a_{1}a_{2}}^{a_{3}} \to 
\mathcal{V}_{a_{1}a'_{3}}^{a'_{2}}$ (see \cite{HL2}).
For $a_{1}, a_{2}, a_{3}\in \A$,
$\mathcal{Y}\in \mathcal{V}_{a_{1}a_{2}}^{a_{3}}$, 
we define 
\begin{eqnarray*}
\sigma_{12}(\mathcal{Y})&=&e^{\pi i \Delta(\mathcal{Y})}
\Omega_{-1}(\mathcal{Y})\\
&=&e^{-\pi i \Delta(\mathcal{Y})}\Omega_{0}(\mathcal{Y}),\\
\sigma_{23}(\mathcal{Y})&=&e^{\pi i h_{a_{1}}}
A_{-1}(\mathcal{Y})\\
&=&e^{-\pi i h_{a_{1}}}A_{0}(\mathcal{Y}),
\end{eqnarray*}
where $\Delta(\Y)=h_{a_{3}}-h_{a_{1}}-h_{a_{2}}$.
By Proposition 1.1 in \cite{H7}, 
they generate an action of $S_{3}$ on 
$\mathcal{V}$.

We now want to choose
a special basis $\mathcal{Y}_{a_{1}a_{2}; i}^{a_{3}}$, $i=1, \dots, 
N_{a_{1}a_{2}}^{a_{3}}$, 
of $\mathcal{V}_{a_{1}a_{2}}^{a_{3}}$ for the triples 
$(a_{1}, a_{2},
a_{3})$ of the forms
$(e, a, a)$, $(a, e, a)$ and $(a, a', e)$ 
where $a\in \A$.
For $a\in \A$, we choose $\Y_{ea; 1}^{a}$ to be the vertex operator 
$Y_{W^{a}}$ defining the module structure on $W^{a}$ and we choose 
$\Y_{ae; 1}^{a}$ to be $\sigma_{12}(\Y_{ea; 1}^{a})$. 
Since
$V'$ as a $V$-module is equivalent to $V$,  we have 
$e'=e$ and we know from Remark 5.3.3 in \cite{FHL} 
that there is a nondegenerate symmetric
invariant 
bilinear form $(\cdot, \cdot)$ on $V$ such that $(\mathbf{1}, 
\mathbf{1})=1$. We identify $V$ and $V'$ using this form. 
We choose $\Y_{aa'; 1}^{e}=\Y_{aa'; 1}^{e'}$
to be the intertwining operator defined using the action of 
$\sigma_{23}$ by
$$\Y_{aa'; 1}^{e'}=\sigma_{23}(\Y_{ae; 1}^{a}),$$
that is, 
$\Y_{aa'; 1}^{e}$ is defined by 
$$(u, \Y_{aa'; 1}^{e}(w_{a}, x)w_{a}')=\langle 
\Y_{ae; 1}^{a}(e^{xL(1)}x^{-2L(0)}w_{a}, x^{-1})u, w_{a}'\rangle$$
for $u\in V$, $w_{a}\in W^{a}$ and $w_{a}'\in (W^{a})'$. 

We now discuss modular transformations. 
Let $q_{\tau}=e^{2\pi i\tau}$
for $\tau \in \mathbb{H}$ ($\mathbb{H}$ is the upper-half plane). 
We consider the $q_{\tau}$-traces of the vertex operators 
$Y_{W^{a}}$ for $a\in \mathcal{A}$ on 
the irreducible $V$-modules $W^{a}$ of the following form:
\begin{equation}\label{1-trace}
\tr_{W^{a}}
Y_{W^{a}}(e^{2\pi iz L(0)}u, e^{2\pi i z})
q_{\tau}^{L(0)-\frac{c}{24}}
\end{equation}
for $u\in V$. 
In \cite{Z}, under some conditions slightly different 
from (mostly stronger than) those we assume in this paper, 
Zhu proved that these $q$-traces are independent of 
$z$, are absolutely convergent
when $0<|q_{\tau}|<1$ and can be analytically extended to 
analytic functions of $\tau$ in the upper-half plane. 
We shall denote the analytic extension of (\ref{1-trace})
by 
$$E(\tr_{W^{a}}
Y_{W^{a}}(e^{2\pi iz L(0)}u, e^{2\pi i z})
q_{\tau}^{L(0)-\frac{c}{24}}).$$
In \cite{Z}, under his conditions alluded to above,
Zhu also proved the following modular invariance property:
For 
$$\left(\begin{array}{cc}a&b\\ c&d\end{array}\right)\in SL(2, \Z),$$
let $\tau'=\frac{a\tau+b}{c\tau+d}$. Then there exist 
unique $A_{a_{1}}^{a_{2}}\in \C$ for $a_{1}, a_{2}\in \mathcal{A}$
such that 
\begin{eqnarray*}
\lefteqn{E\left(\tr_{W^{a_{1}}}
Y_{W^{a_{1}}}\left(e^{\frac{2\pi iz}{c\tau+d}L(0)}
\left(\frac{1}{c\tau+d}\right)^{L(0)}
u, e^{\frac{2\pi i z}{c\tau+d}}\right)
q_{\tau'}^{L(0)-\frac{c}{24}}\right)}\nn
&&=\sum_{a_{2}\in \A}A_{a_{1}}^{a_{2}}
E(\tr_{W^{a_{2}}}
Y_{W^{a_{2}}}(e^{2\pi iz L(0)}u, e^{2\pi i z})
q_{\tau}^{L(0)-\frac{c}{24}})
\end{eqnarray*}
for $u\in V$. In \cite{DLM}, Dong, Li and Mason, 
among many other things,
improved Zhu's results above by showing that the results 
of Zhu above also hold for vertex operator algebras satisfying 
the conditions (slightly weaker than what) we assume in this paper. 
In particular, for 
$$\left(\begin{array}{cc}0&1\\ -1&0\end{array}\right)\in SL(2, \Z),$$
there exist unique $S_{a_{1}}^{a_{2}}\in \C$ for $a_{1}\in \mathcal{A}$
such that 
\begin{eqnarray*}
\lefteqn{E\left(\tr_{W^{a_{1}}}
Y_{W^{a_{1}}}\left(e^{-\frac{2\pi iz}{\tau}L(0)}
\left(-\frac{1}{\tau}\right)^{L(0)}
u, e^{-\frac{2\pi i z}{\tau}}\right)
q_{-\frac{1}{\tau}}^{L(0)-\frac{c}{24}}\right)}\nn
&&=\sum_{a_{2}\in \A}S_{a_{1}}^{a_{2}}
E(\tr_{W^{a_{2}}}
Y_{W^{a_{2}}}(e^{2\pi iz L(0)}u, e^{2\pi i z})
q_{\tau}^{L(0)-\frac{c}{24}})
\end{eqnarray*}
for $u\in V$. When $u=\mathbf{1}$, we see that the matrix
$S=(S_{a_{1}}^{a_{2}})$ actually acts on the space 
of spanned by the 
vacuum characters $\tr_{W^{a}}q_{\tau}^{L(0)-\frac{c}{24}}$
for $a\in \mathcal{A}$.

For a vertex operator algebra $V$ satisfying the conditions 
above, the author proved  in \cite{H7} the Verlinde conjecture 
which states that 
the action of the modular transformation $\tau\mapsto -1/\tau$ 
on the space of characters of irreducible $V$-modules 
diagonalizes the matrices formed by fusion rules. (Note that
in the proof of the Verlinde conjecture, the modular invariance 
results obtained in \cite{Z} and \cite{DLM} are not enough.
One needs the results 
obtained in \cite{H5} and \cite{H6} on intertwining operators and 
on the modular invariance of the space of $q$-traces of products of 
intertwining operators, respectively.)  In this paper,
we need the following two useful consequences 
of the Verlinde conjecture:
For $a\in \A$, 
$$F(\Y_{ae; 1}^{a_{2}} \otimes \Y_{a'a; 1}^{e};
\Y_{ea; 1}^{a}\otimes \Y_{aa'; 1}^{e})\ne 0$$
and 
\begin{equation}\label{s-form-3}
S_{a_{1}}^{a_{2}}
=\frac{S_{e}^{e}(B^{(-1)})_{a_{2}, a_{1}}^{2}}
{F_{a_{1}}
F_{a_{2}}}
\end{equation}
where $(S_{a_{1}}^{a_{2}})$ is the matrix representing 
the action of the modular transformation $\tau\mapsto 
-1/\tau$ on the space of characters of irreducible 
$V$-modules, 
$$(B^{(-1)})_{a_{2}, a_{1}}^{2}
=(B^{(-1)})^{2}(\Y_{a_{2}e; 1}^{a_{2}}
\otimes \Y_{a'_{1}a_{1}; 1}^{e};
\Y_{a_{2}e; 1}^{a_{2}}\otimes \Y_{a'_{1}a_{1}; 1}^{e})
$$
for $a_{1}, a_{2}\in \A$, and 
$$F_{a}=F(\Y_{ae; 1}^{a} \otimes \Y_{a'a; 1}^{e};
\Y_{ea; 1}^{a}\otimes \Y_{aa'; 1}^{e})$$
for $a\in \A$.

\renewcommand{\theequation}{\thesection.\arabic{equation}}
\renewcommand{\thethm}{\thesection.\arabic{thm}}
\setcounter{equation}{0}
\setcounter{thm}{0}

\section{The proof of the rigidity and the dimensions of the
irreducible modules}

In this section, we prove that the braided tensor category structure on
the category of $V$-module is in fact rigid. 

First we recall the definition of rigidity (see for example 
\cite{T} and \cite{BK}).
A tensor category with tensor product bifunctor $\boxtimes$ 
and unit object $V$ 
is rigid if for every object 
$W$ in the category, there are right and left dual objects $W^{*}$
and $^{*}W$ together with morphisms $e_{W}: W^{*}\boxtimes W\to V$,
$i_{W}:V\to W\boxtimes W^{*}$, $e'_{W}: W\boxtimes {}^{*}W\to V$
and $i'_{W}: V\to {}^{*}W\boxtimes W$ 
such that the compositions of the morphisms in the sequences
$$\begin{CD}
W&@>l_{W}^{-1}>>
&V\boxtimes W
&@>i_{W}\boxtimes I_{W}>>
&(W\boxtimes W^{*})\boxtimes 
W&@>>>\\
&@>\A^{-1}>>&W\boxtimes (W^{*}\boxtimes 
W)&@>I_{W}\boxtimes \;e_{W}>>&
W\boxtimes V&@>r_{W}>>&W,
\end{CD}$$
$$\begin{CD}
W@>r^{-1}_{W}>> W\boxtimes V
@>I_{W} \boxtimes \; i'_{W}>>
W\boxtimes ({}^{*}W\boxtimes W)
@>>>  \\
 @>\mathcal{A}>> (W\boxtimes {}^{*}W)\boxtimes 
W@>e'_{W} \boxtimes I_{W}>>
V\boxtimes W @>l_{W}>> W,
\end{CD}$$
$$\begin{CD}
{}^{*}W@>l^{-1}_{{}^{*}\!W}>> V\boxtimes {}^{*}W
@>i'_{W}\boxtimes I_{{}^{*}\!W}>>
({}^{*}W\boxtimes W)\boxtimes 
{}^{*}W@>>>  \\
 @>\mathcal{A}^{-1}>> {}^{*}W\boxtimes (W\boxtimes 
{}^{*}W)@>I_{{}^{*}W}\boxtimes \; e'_{W}>>
{}^{*}W\boxtimes V @>r_{{}^{*}\!W}>> {}^{*}W,
\end{CD}$$
$$\begin{CD}
W^{*}@>r^{-1}_{W^{*}}>> W^{*}\boxtimes V
@>I_{W^{*}} \boxtimes \; i_{W}>>
 W^{*}\boxtimes (W\boxtimes W^{*})
@>>> \\ 
 @>\mathcal{A}>> (W^{*}\boxtimes W)\boxtimes 
W^{*}@>e_{W} \boxtimes I_{W^{*}}>>
V\boxtimes W^{*} @>l_{W^{*}}>> W_{*},
\end{CD}$$
are equal to the identity isomorphisms
$I_{W}$, $I_{W}$, $I_{{}^{*}W}$, $I_{W^{*}}$,
respectively. Rigidity is a standard notion in the theory 
of tensor categories. See, for example, 
\cite{T} and \cite{BK} for details.

In this section, we shall always use the bases
$\{\Y_{ea;1}^{a}\}$, $\{\Y_{ae;1}^{a}\}$ and $\{\Y_{aa';1}^{e}\}$
of $\V_{ea}^{a}$, $\V_{ae}^{a}$ and $\V_{aa'}^{e}$ chosen
in the preceding section.
We take both the left and right duals of a $V$-module $W$
to be the contragredient module $W'$ of $W$. Since our 
tensor category is semisimple, to prove the rigidity, we need 
only discuss irreducible modules. 
For $a\in \mathcal{A}$ and $z\in \C^{\times}$,
using the universal property
for the tensor product module $(W^{a})'\boxtimes_{P(z)} W^{a}$,
we know that there exists a unique module map $\hat{e}_{a; z}:
(W^{a})'\boxtimes_{P(z)} W^{a}\to V$ such that 
$$\overline{\hat{e}_{a; z}}(w'_{a}\boxtimes_{P(z)} w_{a})
=\Y_{a'a; 1}^{e}(w'_{a}, z)w_{a}$$
for $w_{a}\in W^{a}$ and $w'_{a}\in (W^{a})'$, where 
$\overline{\hat{e}_{a; z}}: \overline{(W^{a})'\boxtimes_{P(z)} W^{a}}\to 
\overline{V}$ is the natural extension of $\hat{e}_{a; z}$. 
When $z=1$, we shall denote $\hat{e}_{a; 1}$ simply by $\hat{e}_{a}$.

For  $a\in \mathcal{A}$ and $z\in \C^{\times}$, 
by the universal property 
for tensor product modules, 
for any fixed bases $\{\Y_{aa'; i}^{b}\}$ of $\mathcal{V}_{aa'}^{b}$
for $b\in \mathcal{A}$, there exists an isomorphism 
$$f_{a; z}: W^{a}\boxtimes_{P(z)} (W^{a})'\to 
\coprod_{b\in \mathcal{A}}N_{aa'}^{b}
W^{b}$$
such that 
$$\overline{f_{a; z}}(w_{a}\boxtimes_{P(z)} w'_{a})
=\Y(w_{a}, z)w'_{a}$$
where 
$$\Y=\sum_{b\in \mathcal{A}}\sum_{i=1}^{N_{aa'}^{b}}\Y_{aa'; i}^{b},$$
an intertwining operator of type ${\coprod_{b\in \mathcal{A}}N_{aa'}^{b}
W^{b}\choose W^{a}(W^{a})'}$. Since $N_{aa'}^{e}=1$, there is 
a unique injective module map 
$$g_{a; z}: V \to 
\coprod_{b\in \mathcal{A}}N_{aa'}^{b}
W^{b}$$
such that $g_{a; z}(\mathbf{1})=\mathbf{1}\in V=W^{e}\subset 
\coprod_{b\in \mathcal{A}}N_{aa'}^{b}
W^{b}$. Let 
$$i_{a; z}=f_{a; z}^{-1}\circ 
g_{a; z}: V\to W^{a}\boxtimes_{P(z)} (W^{a})'.$$
When $z=1$, we shall denote $i_{a; 1}$ simply by $i_{a}$.

\begin{lemma}\label{ind-bases}
The map $i_{a; z}$ is independent of the choice of 
the bases $\Y_{aa'; i}^{b}$ of $\mathcal{V}_{aa'}^{b}$
for $b\in \mathcal{A}$, $b\ne e$.
\end{lemma}
\pf
We choose bases $\widetilde{\Y}_{aa'; i}^{b}$ of 
$\mathcal{V}_{aa'}^{b}$
for $b\in \mathcal{A}$ such that 
$\widetilde{\Y}_{aa'; 1}^{e}=\Y_{aa'; 1}^{e}$. 
Then we also obtain an isomorphism
$$\tilde{f}_{a; z}: W^{a}\boxtimes_{P(z)} (W^{a})'\to 
\coprod_{b\in \mathcal{A}}N_{aa'}^{b}
W^{b}.$$
Since $\widetilde{\Y}_{aa'; 1}^{e}=\Y_{aa'; 1}^{e}$, 
$\pi_{a}\circ \tilde{f}_{a; z}=\pi_{a} \circ f_{a; z}$ 
where $\pi_{a}: 
\coprod_{b\in \mathcal{A}}N_{aa'}^{b}
W^{b}\to W^{e}=V$ is the projection.
By the definition of $g_{a; z}$, we see that 
$f_{a; z}^{-1}\circ 
g_{a; z}=\tilde{f}_{a; z}^{-1}\circ 
g_{a; z}$.
So $i_{a; z}$ is independent of the choices of 
$\Y_{aa'; i}^{b}$ of $\mathcal{V}_{aa'}^{b}$
for $b\in \mathcal{A}$, $b\ne e$.
\epfv

For $a\in \mathcal{A}$ and $z\in \C^{\times}$,
using the universal property
for the tensor product module $W^{a} \boxtimes_{P(z)} (W^{a})'$,
we know that there exists a unique module map $\hat{e}'_{a; z}:
W^{a} \boxtimes_{P(z)} (W^{a})'\to V$ such that 
$$\overline{\hat{e}'_{a; z}}(w_{a}\boxtimes_{P(z)} w'_{a})
=\Y_{aa'; 1}^{e}(w_{a}, z)w'_{a}$$
for $w_{a}\in W^{a}$ and $w'_{a}\in (W^{a})'$, where 
$\overline{\hat{e}'_{a; z}}: \overline{W^{a} \boxtimes_{P(z)} (W^{a})'}\to 
\overline{V}$ is the natural extension of $\hat{e}'_{a; z}$. 
When $z=1$, we shall denote $\hat{e}'_{a; 1}$ simply by $\hat{e}'_{a}$.

For any fixed $a\in \mathcal{A}$ and $z\in \C^{\times}$, 
by the universal property 
for tensor product modules, 
for any fixed basis $\Y_{a'a; i}^{b}$ of $\mathcal{V}_{a'a}^{b}$ as we 
choose in Section 1
for $b\in \mathcal{A}$, $b\ne e$, there exists an isomorphism 
$$f'_{a; z}: (W^{a})'\boxtimes_{P(z)} W^{a}\to 
\coprod_{b\in \mathcal{A}}N_{a'a}^{b}
W^{b}$$
such that 
$$\overline{f'_{a; z}}(w'_{a}\boxtimes_{P(z)} w_{a})
=\Y(w'_{a}, z)w_{a}$$
where 
$$\Y=\sum_{b\in \mathcal{A}}\sum_{i=1}^{N_{a'a}^{b}}\Y_{a'a; i}^{b},$$
an intertwining operator of type ${\coprod_{b\in \mathcal{A}}N_{a'a}^{b}
W^{b}\choose (W^{a})'W^{a}}$. Since $N_{a'a}^{e}=1$, there is 
a unique injective module map 
$$g'_{a; z}: V \to 
\coprod_{b\in \mathcal{A}}N_{a'a}^{b}
W^{b}$$
such that $g'_{a; z}(\mathbf{1})=\mathbf{1}\in V=W^{e}\subset 
\coprod_{b\in \mathcal{A}}N_{a'a}^{b}
W^{b}$. Let 
$$i'_{a; z}=(f'_{a; z})^{-1}\circ 
g'_{a; z}: V\to (W^{a})'\boxtimes_{P(z)} W^{a}.$$
When $z=1$, we shall denote $i'_{a; 1}$ simply by $i'_{a}$.

\begin{lemma}
The map $i'_{a; z}$ is independent of the choice of 
the bases $\Y_{a'a; i}^{b}$ of $\mathcal{V}_{a'a}^{b}$
for $b\in \mathcal{A}$, $b\ne e$.
\end{lemma}
\pf
The proof is the same as the one for Lemma 
\ref{ind-bases}.
\epfv

The composition 
$$r_{a}\circ (I_{W^{a}}\boxtimes \hat{e}_{a})\circ 
\mathcal{A}^{-1}\circ (i_{a}\boxtimes I_{W^{a}}) \circ l^{-1}_{a}$$
of the module maps in the sequence
\begin{equation}\label{map-1}
\begin{CD}
W^{a}@>l^{-1}_{a}>> V\boxtimes W^{a}
@>i_{a}\boxtimes I_{W^{a}}>>
(W^{a}\boxtimes (W^{a})')\boxtimes 
W^{a}@>>>  \\
 @>\mathcal{A}^{-1}>> W^{a}\boxtimes ((W^{a})'\boxtimes 
W^{a})@>I_{W^{a}}\boxtimes \; \hat{e}_{a}>>
W^{a}\boxtimes V @>r_{a}>> W_{a}
\end{CD}
\end{equation}
is a module map 
from  $W^{a}$ to itself. Since $W^{a}$ is irreducible,
there must exist $\lambda_{a}\in \C$ such that this module map is 
equal to $\lambda_{a}$ times the identity map on $W^{a}$, that is,
\begin{equation}\label{map1}
r_{a}\circ (I_{W^{a}}\boxtimes \hat{e}_{a})\circ 
\mathcal{A}^{-1}\circ (i_{a}\boxtimes I_{W^{a}}) \circ l^{-1}_{a}
=\lambda_{a}I_{W^{a}}.
\end{equation}
We need to calculate $\lambda_{a}$.

Let $z_{1}, z_{2}$ be any nonzero complex numbers.
We first calculate the composition of 
the module maps given by the sequence
\begin{equation}\label{map'-z}
\begin{CD}
W^{a}@>l_{a; z_{2}}^{-1}>>V\boxtimes_{P(z_{2})} W^{a}\\
@>i_{a; z_{1}-z_{2}}\boxtimes_{P(z_{2})} I_{W^{a}}>>
(W^{a}\boxtimes_{P(z_{1}-z_{2})}
(W^{a})')\boxtimes_{P(z_{2})} W^{a}\\
@>\left(\mathcal{A}_{P(z_{1}), P(z_{2})}^{P(z_{1}-z_{2}), 
P(z_{2})}\right)^{-1}>>
W^{a}\boxtimes_{P(z_{1})} ((W^{a})'\boxtimes_{P(z_{2})} 
W^{a})\\
@>I_{W^{a}}\boxtimes_{P(z_{1})} \; \hat{e}_{a; z_{2}}>>
W^{a}\boxtimes_{P(z_{1})} V@>r_{a; z_{1}}>> W_{a},
\end{CD}
\end{equation}
that is, 
we first calculate
$$r_{a; z_{1}}\circ (I_{W^{a}}\boxtimes_{P(z_{1})}
\hat{e}_{a; z_{2}})\circ
\left(\mathcal{A}_{P(z_{1}), P(z_{2})}^{P(z_{1}-z_{2}),
P(z_{2})}\right)^{-1}\circ
(i_{a; z_{1}-z_{2}}\boxtimes_{P(z_{2})}
I_{W^{a}})\circ l_{a; z_{2}}^{-1}.$$
Since this is also a module map from an irreducible module to 
itself, there exists $\lambda_{a, z_{1}, z_{2}}\in \C$ such that 
it is equal to $\lambda_{a; z_{1}, z_{2}}I_{W^{a}}$.

\begin{prop}\label{main-calc}
For $a\in \mathcal{A}$ and $z_{1}, z_{2}\in \C^{\times}$ satisfying 
$|z_{1}|>|z_{2}|>|z_{1}-z_{2}|>0$, we have
\begin{eqnarray}\label{lambda}
\lambda_{a; z_{1}, z_{2}}&=&F_{a}\nn
&=&F(\Y_{ae; 1}^{a} \otimes \Y_{a'a; 1}^{e};
\Y_{ea; 1}^{a}\otimes \Y_{aa'; 1}^{e}).
\end{eqnarray}
\end{prop}
\pf
Let $w_{1}, w_{2}\in W^{a}$ and $w_{1}'\in (W^{a})'$. 
Then 
$$\Y_{ea;1}^{a}(\Y_{aa';1}^{e}(w_{1}, z_{1}-z_{2})w'_{1}, z_{2})w_{2}\in 
\overline{W^{a}}.$$
By the definition of $l_{a; z_{2}}$, we have
\begin{equation}\label{calc-1}
\overline{l_{a; z_{2}}^{-1}}(
\Y_{ea;1}^{a}(\Y_{aa';1}^{e}(w_{1}, z_{1}-z_{2})w'_{1}, z_{2})w_{2})
=(\Y_{aa';1}^{e}(w_{1}, z_{1}-z_{2})w'_{1})\boxtimes_{P(z_{2})}w_{2}.
\end{equation}

By definitions and the associativity of intertwining operators, we have 
\begin{eqnarray}\label{calc-2}
\lefteqn{\overline{r_{a; z_{1}}}
\biggl(\overline{\left(I_{W^{a}}\boxtimes_{P(z_{1})}
\hat{e}_{a; z_{2}}\right)}
\biggl(\;\overline{\left(\mathcal{A}_{P(z_{1}), P(z_{2})}^{P(z_{1}-z_{2}),
P(z_{2})}\right)^{-1}}
((w_{1}\boxtimes_{P(z_{1}-z_{2})}w_{1}')\boxtimes_{P(z_{2})}
w_{2})\biggr)\biggr)}\nn
&&=\overline{r_{a; z_{1}}}
\left(\overline{\left(I_{W^{a}}\boxtimes_{P(z_{1})}
\hat{e}_{a; z_{2}}\right)}
(w_{1}\boxtimes_{P(z_{1})}(w_{1}'\boxtimes_{P(z_{2})}
w_{2}))\right)\nn
&&=\overline{r_{a; z_{1}}}(w_{1}\boxtimes_{P(z_{1})}
(\Y_{a'a;1}^{e}(w_{1}', z_{2})w_{2}))\nn
&&=\Y_{ae;1}^{a}(w_{1}, z_{1})\Y_{a'a;1}^{e}(w_{1}', z_{2})w_{2}\nn
&&=\sum_{b\in \A}\sum_{k=1}^{N_{ba}^{a}}\sum_{l=1}^{N_{aa'}^{b}}
F(\Y_{ae; 1}^{a} \otimes \Y_{a'a; 1}^{e};
\Y_{ba; k}^{a}\otimes \Y_{aa'; l}^{b})\cdot\nn
&& \quad\quad\quad\quad\quad\quad\quad\quad\quad
\quad\quad\quad\quad\cdot
\Y_{ba;k}^{a}(\Y_{aa';l}^{b}(w_{1}, z_{1}-z_{2})w_{1}', z_{2})w_{2},
\quad\quad\quad\quad
\end{eqnarray}
when $|z_{1}|>|z_{2}|>|z_{1}-z_{2}|>0$.
Since both sides of 
(\ref{calc-2}) are well-defined in the region $|z_{2}|>|z_{1}-z_{2}|>0$,
the left- and right-hand sides of (\ref{calc-2}) are equal 
in this larger region as series in (rational) powers of $z_{2}$
and $z_{1}-z_{2}$. 

Let $\eta_{ba; k}^{a}: W^{b}\boxtimes_{P(z_{2})}W^{a}
\to W^{a}$ and $\eta_{aa';l}^{b}: W^{a}\boxtimes_{P(z_{1}-z_{2})}(W^{a})'
\to W^{b}$ be module maps determined by 
$$\eta_{ba; k}^{a}(w_{3}\boxtimes_{P(z_{2})}w_{2})
=\Y_{ba;k}^{a}(w_{3}, z_{2})w_{2}$$
and 
$$\eta_{aa';l}^{b}(w_{1}\boxtimes_{P(z_{2})}w_{1}')
=\Y_{aa';k}^{a}(w_{1}, z_{1}-z_{2})w_{1}',$$
respectively, for $w_{1}, w_{2}\in W^{a}$, $w_{3}\in W^{b}$ and 
$w_{1}'\in (W^{a})'$. 
Then (\ref{calc-2}) gives
\begin{eqnarray}\label{calc-2.1}
\lefteqn{\overline{r_{a; z_{1}}}
\biggl(\overline{\left(I_{W^{a}}\boxtimes_{P(z_{1})}
\hat{e}_{a; z_{2}}\right)}
\biggl(\;\overline{\left(\mathcal{A}_{P(z_{1}), P(z_{2})}^{P(z_{1}-z_{2}),
P(z_{2})}\right)^{-1}}
((w_{1}\boxtimes_{P(z_{1}-z_{2})}w_{1}')\boxtimes_{P(z_{2})}
w_{2})\biggr)\biggr)}\nn
&&=\sum_{b\in \A}\sum_{k=1}^{N_{ba}^{a}}\sum_{l=1}^{N_{aa'}^{b}}
F(\Y_{ae; 1}^{a} \otimes \Y_{a'a; 1}^{e};
\Y_{ba; k}^{a}\otimes \Y_{aa'; l}^{b})\cdot\nn
&& \quad\quad\cdot
(\overline{\eta_{ba;k}^{a}}\circ (\overline{\eta_{aa';l}^{b}}\boxtimes
_{P(z_{1}-z_{2})} I_{(W^{a})'}))
((w_{1}\boxtimes_{P(z_{1}-z_{2})}w_{1}')\boxtimes_{P(z_{2})}w_{2})
\quad\quad\quad
\end{eqnarray}
for $w_{1}, w_{2}\in W^{a}$ and 
$w_{1}'\in (W^{a})'$, 
when $|z_{2}|>|z_{1}-z_{2}|>0$. Since the components of 
elements of the form 
$(w_{1}\boxtimes_{P(z_{1}-z_{2})}w_{1}')\boxtimes_{P(z_{2})}w_{2}$
for $w_{1}, w_{2}\in W^{a}$ and 
$w_{1}'\in (W^{a})'$ span $(W^{a}\boxtimes_{P(z_{1}-z_{2})}
(W^{a})')\boxtimes_{P(z_{2})}W^{a}$, (\ref{calc-2.1}) gives
\begin{eqnarray*}
\lefteqn{r_{a; z_{1}}\circ 
\left(I_{W^{a}}\boxtimes_{P(z_{1})}
\hat{e}_{a; z_{2}}\right)\circ
\left(\mathcal{A}_{P(z_{1}), P(z_{2})}^{P(z_{1}-z_{2}),
P(z_{2})}\right)^{-1}}\nn
&&=\sum_{b\in \A}\sum_{k=1}^{N_{ba}^{a}}\sum_{l=1}^{N_{aa'}^{b}}
F(\Y_{ae; 1}^{a} \otimes \Y_{a'a; 1}^{e};
\Y_{ba; k}^{a}\otimes \Y_{aa'; l}^{b})\cdot\nn
&& \quad\quad\quad\quad\quad\cdot
(\eta_{ba;k}^{a}\circ (\eta_{aa';l}^{b}\boxtimes
_{P(z_{1}-z_{2})} I_{(W^{a})'}))
\end{eqnarray*}
Thus we obtain
\begin{eqnarray}\label{calc-2.2}
\lefteqn{r_{a; z_{1}}\circ 
\left(I_{W^{a}}\boxtimes_{P(z_{1})}
\hat{e}_{a; z_{2}}\right)\circ
\left(\mathcal{A}_{P(z_{1}), P(z_{2})}^{P(z_{1}-z_{2}),
P(z_{2})}\right)^{-1}\circ (i_{a; z_{1}-z_{2}}\boxtimes_{P(z_{2})}
I_{(W^{a})'})\circ l_{a; z_{2}}^{-1}}\nn
&&=\sum_{b\in \A}\sum_{k=1}^{N_{ba}^{a}}\sum_{l=1}^{N_{aa'}^{b}}
F(\Y_{ae; 1}^{a} \otimes \Y_{a'a; 1}^{e};
\Y_{ba; k}^{a}\otimes \Y_{aa'; l}^{b})\cdot\nn
&& \quad\quad\quad\quad\cdot
(\eta_{ba;k}^{a}\circ (\eta_{aa';l}^{b}\boxtimes
_{P(z_{1}-z_{2})} I_{(W^{a})'})\circ 
(i_{a; z_{1}-z_{2}}\boxtimes_{P(z_{2})}
I_{(W^{a})'})\circ l_{a; z_{2}}^{-1})\nn
&&=\sum_{b\in \A}\sum_{k=1}^{N_{ba}^{a}}\sum_{l=1}^{N_{aa'}^{b}}
F(\Y_{ae; 1}^{a} \otimes \Y_{a'a; 1}^{e};
\Y_{ba; k}^{a}\otimes \Y_{aa'; l}^{b})\cdot\nn
&& \quad\quad\quad\quad\cdot
(\eta_{ba;k}^{a}\circ ((\eta_{aa';l}^{b}\circ i_{a; z_{1}-z_{2}})\boxtimes
_{P(z_{1}-z_{2})} I_{(W^{a})'})\circ  l_{a; z_{2}}^{-1}).
\end{eqnarray}

Let $\pi_{b; l}: \coprod_{d\in \A}N_{aa'}^{d}W^{d}\to 
W^{b}$ for $b\in \A$ and $l=1, 
\dots, N_{aa'}^{b}$ be the projection to the $l$-th copy of $W^{b}$.
Then
\begin{eqnarray*}
\eta_{aa';l}^{b}(w_{1}\boxtimes_{P(z_{1}-z_{2})}w'_{1})
&=&\Y_{aa';l}^{b}(w_{1}, z_{1}-z_{2})w'_{1}\nn
&=&(\pi_{b; l}\circ f_{a; z_{1}-z_{2}})
(w_{1}\boxtimes_{P(z_{1}-z_{2})}w'_{1})
\end{eqnarray*}
for $w_{1}\in W^{a}$ and $w'_{1}\in (W^{a})'$.
So we have 
$\eta_{aa';l}^{b}=\pi_{b; l}\circ f_{a; z_{1}-z_{2}}$ and
\begin{eqnarray}\label{calc-2.3}
\eta_{aa';l}^{b}\circ i_{a; z_{1}-z_{2}}
&=&\pi_{b; l}\circ f_{a; z_{1}-z_{2}}
\circ f^{-1}_{a; z_{1}-z_{2}}\circ g_{a; z_{1}-z_{2}}\nn
&=&\pi_{b; l}\circ g_{a; z_{1}-z_{2}}\nn
&=&\delta_{eb}\delta_{1l}I_{V}.
\end{eqnarray}
Also,
\begin{eqnarray*}
(\eta_{ea;1}^{a}\circ l_{a; z_{2}}^{-1})(w_{1})
&=&\eta_{ea;1}^{a}(l_{a; z_{2}}^{-1}(w_{1}))\nn
&=&\eta_{ea;1}^{a}(\mathbf{1}\boxtimes_{P(z_{2})} w_{1})\nn
&=&Y_{W^{a}}(\mathbf{1}, z_{2})w_{1}\nn
&=&w_{1}
\end{eqnarray*}
for $w_{1}\in W^{a}$. So we have 
\begin{equation}\label{calc-2.4}
\eta_{ea;1}^{a}\circ l_{a; z_{2}}^{-1}
=I_{W^{a}}.
\end{equation}
Using (\ref{calc-2.3}) and (\ref{calc-2.4}),
the right-hand side of (\ref{calc-2.2}) becomes
\begin{eqnarray}\label{calc-2.5}
\lefteqn{\sum_{b\in \A}\sum_{k=1}^{N_{ba}^{a}}\sum_{l=1}^{N_{aa'}^{b}}
F(\Y_{ae; 1}^{a} \otimes \Y_{a'a; 1}^{e};
\Y_{ba; k}^{a}\otimes \Y_{aa'; l}^{b})\cdot}\nn
&& \quad\quad\quad\quad\cdot \delta_{eb}\delta_{1l}
(\eta_{ba;k}^{a}\circ  l_{a; z_{2}}^{-1})\nn
&&=F(\Y_{ae; 1}^{a} \otimes \Y_{a'a; 1}^{e};
\Y_{ea; 1}^{a}\otimes \Y_{aa'; 1}^{e})\cdot\nn
&& \quad\quad\quad\quad\cdot 
(\eta_{ea;1}^{a}\circ  l_{a; z_{2}}^{-1})\nn
&&=F(\Y_{ae; 1}^{a} \otimes \Y_{a'a; 1}^{e};
\Y_{ea; 1}^{a}\otimes \Y_{aa'; 1}^{e})I_{W^{a}}.
\end{eqnarray}

From (\ref{calc-2.2}), (\ref{calc-2.5}) and the definition of 
$\lambda_{a; z_{1}, z_{2}}$, we obtain (\ref{lambda}).
\epfv

\begin{prop}\label{main-cor}
For $a\in \mathcal{A}$, we have
\begin{eqnarray}\label{lambda-1}
\lambda_{a}&=&F_{a}\nn
&=&F(\Y_{ae; 1}^{a} \otimes \Y_{a'a; 1}^{e};
\Y_{ea; 1}^{a}\otimes \Y_{aa'; 1}^{e}).
\end{eqnarray}
\end{prop}
\pf
Choose 
$z_{1}^{0}, z_{2}^{0}\in (0, \infty)$ satisfying
$z_{1}^{0}>z_{2}^{0}>z_{1}^{0}-z_{2}^{0}>0$
and choose paths $\gamma_{0}, \gamma_{1}, \gamma_{2}$ 
in $(0, \infty)$
from $z_{1}^{0}-z_{2}^{0}, z_{1}^{0}, z_{2}^{0}$ to $1, 1, 1$,
respectively.
Then by the definitions of all the module maps involved, 
the following diagrams are commutative:
$$\begin{CD}
W^{a}@>l_{a; z_{2}^{0}}^{-1}>>V\boxtimes_{P(z_{2}^{0})} W^{a}\\
@V=VV @VV\mathcal{T}_{\gamma_{2}}V\\
W^{a}@>l_{a}^{-1}>>V\boxtimes W^{a}
\end{CD}
$$
$$
\begin{CD}
V\boxtimes_{P(z_{2})} W^{a}
@>i_{a; z_{1}^{0}-z_{2}^{0}}\boxtimes_{P(z_{2}^{0})} I_{W^{a}}>>
(W^{a}\boxtimes_{P(z_{1}^{0}-z_{2}^{0})}
(W^{a})')\boxtimes_{P(z_{2}^{0})} W^{a}\\
@V\mathcal{T}_{\gamma_{2}}VV @VV
(\mathcal{T}_{\gamma_{0}}\boxtimes I_{W^{a}})\circ 
\mathcal{T}_{\gamma_{2}}V\\
V\boxtimes W^{a}
@>i_{a}\boxtimes I_{W^{a}}>>
(W^{a}\boxtimes (W^{a})')\boxtimes W^{a}
\end{CD}
$$
$$
\begin{CD}
(W^{a}\boxtimes_{P(z_{1}^{0}-z_{2}^{0})}
(W^{a})')\boxtimes_{P(z_{2}^{0}} W^{a}
@>\left(\mathcal{A}_{P(z_{1}^{0}), P(z_{2}^{0})}^{P(z_{1}^{0}-z_{2}^{0}), 
P(z_{2}^{0})}\right)^{-1}>>
W^{a}\boxtimes_{P(z_{1}^{0})} ((W^{a})'\boxtimes_{P(z_{2}^{0})} 
W^{a})\\
@V(\mathcal{T}_{\gamma_{0}}\boxtimes I_{W^{a}})\circ 
\mathcal{T}_{\gamma_{2}}VV 
@VV(I_{W^{a}}\boxtimes \mathcal{T}_{\gamma_{2}})\circ 
\mathcal{T}_{\gamma_{1}}V\\
(W^{a}\boxtimes
(W^{a})')\boxtimes W^{a}
@>\mathcal{A}^{-1}>>
W^{a}\boxtimes ((W^{a})'\boxtimes
W^{a})
\end{CD}
$$
$$
\begin{CD}
W^{a}\boxtimes_{P(z_{1}^{0})} ((W^{a})'\boxtimes_{P(z_{2}^{0})} 
W^{a})@>I_{W^{a}}\boxtimes_{P(z_{1}^{0})} \hat{e}_{a; z_{2}^{0}}>>
W^{a}\boxtimes_{P(z_{1}^{0})} V\\
@V(I_{W^{a}}\boxtimes \mathcal{T}_{\gamma_{2}})\circ 
\mathcal{T}_{\gamma_{1}}VV
@VV\mathcal{T}_{\gamma_{1}}V\\
W^{a}\boxtimes ((W^{a})'\boxtimes
W^{a})@>I_{W^{a}}\boxtimes \; \hat{e}_{a}>>
W^{a}\boxtimes V
\end{CD}
$$
$$
\begin{CD}
W^{a}\boxtimes_{P(z_{1}^{0})} V@>r_{a; z_{1}^{0}}>>
W^{a}\\
@VV\mathcal{T}_{\gamma_{1}}V @VV=V\\
W^{a}\boxtimes V@>r_{a}>> W^{a}
\end{CD}
$$
Combining these diagrams, we obtain
\begin{eqnarray*}
\lefteqn{r_{a}\circ (I_{W^{a}}\boxtimes \hat{e}_{a})\circ 
\mathcal{A}^{-1}\circ (i_{a}\boxtimes I_{W^{a}}) \circ l^{-1}_{a}}\nn
&&=r_{a; z_{1}}\circ (I_{W^{a}}\boxtimes_{P(z_{1})}
\hat{e}_{a; z_{2}})\;\circ\nn
&&\quad\quad\quad\quad\quad\quad\circ
\left(\mathcal{A}_{P(z_{1}), P(z_{2})}^{P(z_{1}-z_{2}),
P(z_{2})}\right)^{-1}\circ
(i_{a; z_{1}-z_{2}}\boxtimes_{P(z_{2})}
I_{W^{a}})\circ l_{a; z_{2}}^{-1}.
\end{eqnarray*} 
From this equality and the definitions of $\lambda_{a}$ and 
$\lambda_{a; z_{1}^{0}, z_{2}^{0}}$,
we obtain
$$\lambda_{a}=\lambda_{a; z_{1}^{0}, z_{2}^{0}}.$$
By Proposition \ref{main-calc}, we obtain (\ref{lambda-1}).
\epfv

The composition 
$$l_{a}\circ (\hat{e}'_{a} \boxtimes I_{W^{a}})\circ 
\mathcal{A}\circ (I_{W^{a}} \boxtimes i'_{a}) \circ r^{-1}_{a}$$
of the module maps in the sequence
\begin{equation}\label{map-2}
\begin{CD}
W^{a}@>r^{-1}_{a}>> W^{a}\boxtimes V
@>I_{W^{a}} \boxtimes \; i'_{a}>>
W^{a}\boxtimes ((W^{a})'\boxtimes W^{a})
@>>>  \\
 @>\mathcal{A}>> (W^{a}\boxtimes (W^{a})')\boxtimes 
W^{a}@>\hat{e}'_{a} \boxtimes I_{W^{a}}>>
V\boxtimes W^{a} @>l_{a}>> W^{a}
\end{CD}
\end{equation}
is a module map 
from  $W^{a}$ to itself. Again, since $W^{a}$ is irreducible,
there must exist $\mu_{a}\in \C$ such that this module map is 
equal to $\mu_{a}$ times the identity map on $W^{a}$, that is,
\begin{equation}\label{map2}
l_{a}\circ (\hat{e}'_{a} \boxtimes I_{W^{a}})\circ 
\mathcal{A}\circ (I_{W^{a}} \boxtimes i'_{a}) \circ r^{-1}_{a}
=\mu_{a}I_{W^{a}}.
\end{equation}
We need to calculate $\mu_{a}$.

Our method to calculate $\mu_{a}$ is similar to the 
one used to calculate $\lambda_{a}$ above. 
Let $z_{1}, z_{2}$ be any nonzero complex numbers.
We first calculate the composition of 
the module maps given by the sequence
\begin{equation}\label{map'-z-2}
\begin{CD}
W^{a}@>r_{a; z_{2}}^{-1}>>W^{a}\boxtimes_{P(z_{1})} V\\
@>I_{W^{a}}\boxtimes_{P(z_{1})} \;i'_{a; z_{2}}>>
(W^{a}\boxtimes_{P(z_{1})}
((W^{a})'\boxtimes_{P(z_{2})} W^{a})\\
@>\mathcal{A}_{P(z_{1}), P(z_{2})}^{P(z_{1}-z_{2}), 
P(z_{2})}>>
(W^{a}\boxtimes_{P(z_{1}-z_{2})} (W^{a})')\boxtimes_{P(z_{2})} 
W^{a}\\
@>\hat{e}'_{a; z_{1}-z_{2}}\boxtimes_{P(z_{1})} I_{W^{a}}>>
V\boxtimes_{P(z_{2})} W^{a}@>l_{a; z_{2}}>> W^{a},
\end{CD}
\end{equation}
that is, 
we first calculate
$$l_{a; z_{2}}\circ (\hat{e}'_{a; z_{1}-z_{2}}\boxtimes_{P(z_{1})} 
I_{W^{a}})\circ
\mathcal{A}_{P(z_{1}), P(z_{2})}^{P(z_{1}-z_{2}),
P(z_{2})}\circ
(I_{W^{a}}\boxtimes_{P(z_{1})} i'_{a; z_{2}})\circ r_{a; z_{2}}^{-1}.$$
Since this is also a module map from an irreducible module to 
itself, there exists $\mu_{a, z_{1}, z_{2}}\in \C$ such that 
it is equal to $\mu_{a; z_{1}, z_{2}}I_{W^{a}}$.

\begin{prop}\label{main-calc-2}
For $a\in \mathcal{A}$ and $z_{1}, z_{2}\in \C^{\times}$ satisfying 
$|z_{1}|>|z_{2}|>|z_{1}-z_{2}|>0$, we have
\begin{eqnarray}\label{mu}
\mu_{a; z_{1}, z_{2}}&=&F_{a}\nn
&=&F(\Y_{ae; 1}^{a} \otimes \Y_{a'a; 1}^{e};
\Y_{ea; 1}^{a}\otimes \Y_{aa'; 1}^{e}).
\end{eqnarray}
\end{prop}
\pf
Let $w_{1}, w_{2}\in W^{a}$ and $w_{1}'\in (W^{a})'$. 
Then 
$$\Y_{ae;1}^{a}(w_{1}, z_{1})\Y_{a'a;1}^{e}(w'_{1}, z_{2})w_{2}\in 
\overline{W^{a}}.$$
By the definition of $r_{a; z_{2}}$, we have
\begin{equation}\label{calc-2-1}
\overline{r_{a; z_{2}}^{-1}}(
\Y_{ae;1}^{a}(w_{1}, z_{1})\Y_{a'a;1}^{e}(w'_{1}, z_{2})w_{2})
=w_{1}\boxtimes_{P(z_{1})}\Y_{a'a;1}^{e}(w'_{1}, z_{2})w_{2}).
\end{equation}

By definitions and the associativity of intertwining operators, we have 
\begin{eqnarray}\label{calc-2-2}
\lefteqn{\overline{l_{a; z_{2}}}
\biggl(\overline{\left(\hat{e}'_{a; z_{1}-z_{2}}\boxtimes_{P(z_{1})} 
I_{W^{a}}\right)}
\biggl(\;\overline{\mathcal{A}_{P(z_{1}), P(z_{2})}^{P(z_{1}-z_{2}),
P(z_{2})}}
(w_{1}\boxtimes_{P(z_{1})}(w_{1}'\boxtimes_{P(z_{2})}
w_{2}))\biggr)\biggr)}\nn
&&=\overline{l_{a; z_{2}}}
\left(\overline{\left(\hat{e}'_{a; z_{1}-z_{2}}\boxtimes_{P(z_{1})} 
I_{W^{a}}\right)}
((w_{1}\boxtimes_{P(z_{1}-z_{2})}w_{1}')\boxtimes_{P(z_{2})}
w_{2})\right)\nn
&&=\overline{l_{a; z_{2}}}((\Y_{aa';1}^{e}(w_{1}, z_{1}-z_{2})
w_{1}')\boxtimes_{P(z_{2})}w_{2})\nn
&&=\Y_{ea;1}^{a}(\Y_{aa';1}^{e}(w_{1}, z_{1}-z_{2})w_{1}', z_{2})w_{2}\nn
&&=\sum_{b\in \A}\sum_{k=1}^{N_{ab}^{a}}\sum_{l=1}^{N_{a'a}^{b}}
F^{-1}(\Y_{ea; 1}^{a}\otimes \Y_{aa'; 1}^{e};
\Y_{ab; k}^{a} \otimes \Y_{a'a; l}^{b})\cdot\nn
&& \quad\quad\quad\quad\quad\quad\quad\quad\quad
\quad\quad\quad\quad\cdot
\Y_{ab;k}^{a}(w_{1}, z_{1})\Y_{a'a;l}^{b}(w'_{1}, z_{2})w_{2},\quad\quad
\end{eqnarray}
when $|z_{1}|>|z_{2}|>|z_{1}-z_{2}|>0$.
Since both sides of 
(\ref{calc-2-2}) are well-defined in the region $|z_{1}|>|z_{2}|>0$,
the left- and right-hand sides of (\ref{calc-2-2}) are equal 
in this larger region as series in (rational) powers of 
$z_{1}$ and $z_{2}$. 

Let $\eta_{ab; k}^{a}: W^{a}\boxtimes_{P(z_{1})}W^{b}
\to W^{a}$ and $\eta_{a'a;l}^{b}: (W^{a})'\boxtimes_{P(z_{2})}W^{a}
\to W^{b}$ be module maps determined by 
$$\eta_{ab; k}^{a}(w_{1}\boxtimes_{P(z_{1})}w_{3})
=\Y_{ab;k}^{a}(w_{1}, z_{2})w_{3}$$
and 
$$\eta_{a'a;l}^{b}(w'_{1}\boxtimes_{P(z_{2})}w_{2})
=\Y_{a'a;k}^{b}(w'_{1}, z_{2})w_{2},$$
respectively, for $w_{1}, w_{2}\in W^{a}$, $w_{3}\in W^{b}$ and 
$w_{1}'\in (W^{a})'$. 
Then (\ref{calc-2-2}) gives
\begin{eqnarray}\label{calc-2-2.1}
\lefteqn{\overline{l_{a; z_{2}}}
\biggl(\overline{\left(\hat{e}'_{a; z_{1}-z_{2}}\boxtimes_{P(z_{1})} 
I_{W^{a}}\right)}
\biggl(\;\overline{\mathcal{A}_{P(z_{1}), P(z_{2})}^{P(z_{1}-z_{2}),
P(z_{2})}}
(w_{1}\boxtimes_{P(z_{1})}(w_{1}'\boxtimes_{P(z_{2})}
w_{2}))\biggr)\biggr)}\nn
&&=\sum_{b\in \A}\sum_{k=1}^{N_{ab}^{a}}\sum_{l=1}^{N_{a'a}^{b}}
F^{-1}(\Y_{ea; 1}^{a}\otimes \Y_{aa'; 1}^{e};
\Y_{ab; k}^{a} \otimes \Y_{a'a; l}^{b})\cdot\nn
&& \quad\quad\cdot
(\overline{\eta_{ab;k}^{a}}\circ (\overline{I_{W^{a}}\boxtimes
_{P(z_{1})} \eta_{a'a;l}^{b}} ))
(w_{1}\boxtimes_{P(z_{1})}(w_{1}'\boxtimes_{P(z_{2})}w_{2}))
\quad\quad\quad
\end{eqnarray}
for $w_{1}, w_{2}\in W^{a}$ and 
$w_{1}'\in (W^{a})'$, 
when $|z_{1}|>|z_{2}|>0$. Since the components of 
elements of the form 
$w_{1}\boxtimes_{P(z_{1})}(w_{1}'\boxtimes_{P(z_{2})}w_{2})$
for $w_{1}, w_{2}\in W^{a}$ and 
$w_{1}'\in (W^{a})'$ span $W^{a}\boxtimes_{P(z_{1})}
((W^{a})'\boxtimes_{P(z_{2})}W^{a})$, (\ref{calc-2-2.1}) gives
\begin{eqnarray*}
\lefteqn{l_{a; z_{2}}\circ
\left(\hat{e}'_{a; z_{1}-z_{2}}\boxtimes_{P(z_{1})} 
I_{W^{a}}\right)\circ
\mathcal{A}_{P(z_{1}), P(z_{2})}^{P(z_{1}-z_{2}),
P(z_{2})}}\nn
&&=\sum_{b\in \A}\sum_{k=1}^{N_{ab}^{a}}\sum_{l=1}^{N_{a'a}^{b}}
F^{-1}(\Y_{ea; 1}^{a}\otimes \Y_{aa'; 1}^{e};
\Y_{ab; k}^{a} \otimes \Y_{a'a; l}^{b})\cdot\nn
&& \quad\quad\quad\quad\quad\cdot
(\eta_{ab;k}^{a}\circ (I_{W^{a}}\boxtimes
_{P(z_{1})} \eta_{a'a;l}^{b}))
\end{eqnarray*}
Thus we obtain
\begin{eqnarray}\label{calc-2-2.2}
\lefteqn{l_{a; z_{2}}\circ
\left(\hat{e}'_{a; z_{1}-z_{2}}\boxtimes_{P(z_{1})} 
I_{W^{a}}\right)\circ
\mathcal{A}_{P(z_{1}), P(z_{2})}^{P(z_{1}-z_{2}),
P(z_{2})}\circ
(I_{W^{a}}\boxtimes_{P(z_{1})} i'_{a; z_{2}})
\circ r_{a; z_{2}}^{-1}}\nn
&&=\sum_{b\in \A}\sum_{k=1}^{N_{ab}^{a}}\sum_{l=1}^{N_{a'a}^{b}}
F^{-1}(\Y_{ea; 1}^{a}\otimes \Y_{aa'; 1}^{e};
\Y_{ab; k}^{a} \otimes \Y_{a'a; l}^{b})\cdot\nn
&& \quad\quad\quad\quad\quad\cdot
(\eta_{ab;k}^{a}\circ (I_{W^{a}}\boxtimes
_{P(z_{1})} \eta_{a'a;l}^{b})\circ
(I_{W^{a}}\boxtimes_{P(z_{1})} i'_{a; z_{2}})
\circ r_{a; z_{2}}^{-1})\nn
&&=\sum_{b\in \A}\sum_{k=1}^{N_{ab}^{a}}\sum_{l=1}^{N_{a'a}^{b}}
F^{-1}(\Y_{ea; 1}^{a}\otimes \Y_{aa'; 1}^{e};
\Y_{ab; k}^{a} \otimes \Y_{a'a; l}^{b})\cdot\nn
&& \quad\quad\quad\quad\quad\cdot
(\eta_{ab;k}^{a}\circ (I_{W^{a}}\boxtimes
_{P(z_{1})} (\eta_{a'a;l}^{b}\circ i'_{a; z_{2}}))
\circ r_{a; z_{2}}^{-1}).
\end{eqnarray}

Let $\pi'_{b; l}: \coprod_{d\in \A}N_{a'a}^{d}W^{d}\to 
W^{b}$ for $b\in \A$ and $l=1, 
\dots, N_{a'a}^{b}$ be the projection to the $l$-th copy of $W^{b}$.
Then from
\begin{eqnarray*}
\eta_{a'a;l}^{b}(w'_{1}\boxtimes_{P(z_{2})}w_{1})
&=&\Y_{aa';l}^{b}(w'_{1}, z_{2})w_{1}\nn
&=&(\pi_{b; l}\circ f'_{a; z_{2}})
(w'_{1}\boxtimes_{P(z_{1}-z_{2})}w_{1})
\end{eqnarray*}
for $w_{1}\in W^{a}$ and $w'_{1}\in (W^{a})'$, we obtain
$\eta_{a'a;l}^{b}=\pi_{b; l}\circ f'_{a; z_{2}}$ and
\begin{eqnarray}\label{calc-2-2.3}
\eta_{a'a;l}^{b}\circ i'_{a; z_{2}}
&=&\pi_{b; l}\circ f'_{a; z_{2}}
\circ (f'_{a; z_{2}})^{-1}\circ g'_{a; z_{2}}\nn
&=&\pi_{b; l}\circ g'_{a; z_{2}}\nn
&=&\delta_{eb}\delta_{1l}I_{V}.
\end{eqnarray}
We also have
\begin{eqnarray}\label{calc-2-2.4}
(\overline{\eta_{ae;1}^{a}\circ r_{a; z_{2}}^{-1}})(e^{z_{2}L(-1)}w_{1})
&=&\overline{\eta_{ae;1}^{a}}
(\overline{r_{a; z_{2}}^{-1}}(e^{z_{2}L(-1)}w_{1}))\nn
&=&\overline{\eta_{ae;1}^{a}}(w_{1}\boxtimes_{P(z_{2})} \mathbf{1})\nn
&=&\Y_{ae;1}^{a}(w_{1}, z_{2})\mathbf{1}\nn
&=&e^{z_{2}L(-1)}Y_{W^{a}}(\mathbf{1}, -z_{2})w_{1}\nn
&=&e^{z_{2}L(-1)}w_{1}
\end{eqnarray}
for $w_{1}\in W^{a}$. Since the homogeneous components of 
elements of the form $e^{z_{2}L(-1)}w_{1}$ span $W^{a}$, 
we see that (\ref{calc-2-2.4}) gives
\begin{equation}\label{calc-2-2.5}
\eta_{ae;1}^{a}\circ r_{a; z_{2}}^{-1}
=I_{W^{a}}.
\end{equation}
From (\ref{calc-2-2.3}) and (\ref{calc-2-2.5}),
the right-hand side of (\ref{calc-2-2.2}) becomes
\begin{eqnarray}\label{calc-2-2.6}
\lefteqn{\sum_{b\in \A}\sum_{k=1}^{N_{ab}^{a}}\sum_{l=1}^{N_{a'a}^{b}}
F^{-1}(\Y_{ea; 1}^{a}\otimes \Y_{aa'; 1}^{e};
\Y_{ab; k}^{a} \otimes \Y_{a'a; l}^{b})\cdot}\nn
&& \quad\quad\quad\quad\quad\cdot \delta_{eb}\delta_{1l}
(\eta_{ab;k}^{a}\circ  r_{a; z_{2}}^{-1})\nn
&&=F^{-1}(\Y_{ea; 1}^{a}\otimes \Y_{aa'; 1}^{e};
\Y_{ab; k}^{a} \otimes \Y_{a'a; l}^{b})
(\eta_{ae;1}^{a}\circ  r_{a; z_{2}}^{-1})\nn
&&=F^{-1}(\Y_{ea; 1}^{a}\otimes \Y_{aa'; 1}^{e};
\Y_{ab; k}^{a} \otimes \Y_{a'a; l}^{b})I_{W^{a}}.
\end{eqnarray}

From (\ref{calc-2-2.2}), (\ref{calc-2-2.6}) and the definition of 
$\mu_{a; z_{1}, z_{2}}$, we obtain 
\begin{equation}\label{calc-2-4}
\mu_{a; z_{1}, z_{2}}=F^{-1}(\Y_{ae; 1}^{a} \otimes \Y_{a'a; 1}^{e};
\Y_{ea; 1}^{a}\otimes \Y_{aa'; 1}^{e}).
\end{equation}
By Proposition 3.2  in \cite{H7} and the properties of 
the bases $\Y_{ae; 1}^{a}$, $\Y_{a'a; 1}^{e}$, $\Y_{ea; 1}^{a}$
and $\Y_{aa'; 1}^{e}$ under the action of the $S_{3}$ symmetry 
for $a\in \A$, we have 
\begin{equation}\label{calc-2-5}
F^{-1}(\Y_{ae; 1}^{a} \otimes \Y_{a'a; 1}^{e};
\Y_{ea; 1}^{a}\otimes \Y_{aa'; 1}^{e})
=F(\Y_{ae; 1}^{a} \otimes \Y_{a'a; 1}^{e};
\Y_{ea; 1}^{a}\otimes \Y_{aa'; 1}^{e}).
\end{equation}
From (\ref{calc-2-4}) and (\ref{calc-2-5}), we obtain (\ref{mu}).
\epfv

\begin{prop}\label{main-cor-2}
For $a\in \mathcal{A}$, we have
\begin{eqnarray}\label{mu-1}
\mu_{a}&=&F_{a}\nn
&=&F(\Y_{ae; 1}^{a} \otimes \Y_{a'a; 1}^{e};
\Y_{ea; 1}^{a}\otimes \Y_{aa'; 1}^{e}).
\end{eqnarray}
\end{prop}
\pf
Choose 
$z_{1}^{0}, z_{2}^{0}\in (0, \infty)$ satisfying
$z_{1}^{0}>z_{2}^{0}>z_{1}^{0}-z_{2}^{0}>0$
and choose paths $\gamma_{0}, \gamma_{1}, \gamma_{2}$ 
in $(0, \infty)$
from $z_{1}^{0}-z_{2}^{0}, z_{1}^{0}, z_{2}^{0}$ to $1, 1, 1$,
respectively.
Then by the definitions of all the module maps involved for 
general $z_{1}, z_{2}$, the following diagrams are commutative:
$$\begin{CD}
W^{a}@>r_{a; z_{2}^{0}}^{-1}>>W^{a}\boxtimes_{P(z_{2}^{0})}V \\
@V=VV @VV\mathcal{T}_{\gamma_{2}}V\\
W^{a}@>r_{a}^{-1}>>W^{a}\boxtimes V
\end{CD}
$$
$$
\begin{CD}
W^{a}\boxtimes_{P(z^{0}_{2})} V
@>I_{W^{a}}\boxtimes_{P(z^{0}_{1})} \; i'_{a; z_{2}}>>
W^{a}\boxtimes_{P(z_{1}^{0})}
((W^{a})'\boxtimes_{P(z_{2}^{0})} W^{a})\\
@V\mathcal{T}_{\gamma_{2}}VV @VV
(I_{W^{a}}\boxtimes
\mathcal{T}_{\gamma_{2}}) \circ \mathcal{T}_{\gamma_{1}}V\\
W^{a}\boxtimes V
@>I_{W^{a}} \boxtimes \; i'_{a} >>
W^{a}\boxtimes ((W^{a})'\boxtimes W^{a})
\end{CD}
$$
$$
\begin{CD}
W^{a}\boxtimes_{P(z_{1}^{0})}
((W^{a})'\boxtimes_{P(z_{2}^{0})} W^{a})
@>\mathcal{A}_{P(z_{1}^{0}), P(z_{2}^{0})}^{P(z_{1}^{0}-z_{2}^{0}), 
P(z_{2}^{0})}>>
(W^{a}\boxtimes_{P(z_{1}^{0}-z_{2}^{0})} (W^{a})')\boxtimes_{P(z_{2}^{0})} 
W^{a}\\
@V(I_{W^{a}}\boxtimes
\mathcal{T}_{\gamma_{2}}) \circ \mathcal{T}_{\gamma_{1}}VV 
@VV(\mathcal{T}_{\gamma_{0}}\boxtimes I_{W^{a}})\circ 
\mathcal{T}_{\gamma_{2}}V\\
W^{a}\boxtimes ((W^{a})'\boxtimes W^{a})
@>\mathcal{A}>>
(W^{a}\boxtimes (W^{a})')\boxtimes
W^{a}
\end{CD}
$$
$$
\begin{CD}
(W^{a}\boxtimes_{P(z_{1}^{0}-z_{2}^{0})} (W^{a})')\boxtimes_{P(z_{2}^{0})} 
W^{a}
@>\hat{e}'_{a; z_{1}-z_{2}}\boxtimes_{P(z_{1})} 
I_{W^{a}}>>
V\boxtimes_{P(z_{2}^{0})} W^{a}\\
@V(\mathcal{T}_{\gamma_{0}}\boxtimes I_{W^{a}})\circ 
\mathcal{T}_{\gamma_{2}}VV
@VV\mathcal{T}_{\gamma_{2}}V\\
(W^{a}\boxtimes (W^{a})')\boxtimes
W^{a}@>\hat{e}'_{a}\boxtimes I_{W^{a}} >>
V\boxtimes W^{a}
\end{CD}
$$
$$
\begin{CD}
V\boxtimes_{P(z_{2}^{0})} W^{a}@>l_{a; z_{2}^{0}}>>
W^{a}\\
@V\mathcal{T}_{\gamma_{2}}VV @VV=V\\
V\boxtimes W^{a}@>l_{a}>> W^{a}
\end{CD}
$$
Combining these diagrams, we obtain
\begin{eqnarray*}
\lefteqn{l_{a}\circ (\hat{e}'_{a} \boxtimes I_{W^{a}})\circ 
\mathcal{A}\circ (I_{W^{a}} \boxtimes i'_{a}) \circ r^{-1}_{a}}\nn
&&=l_{a; z_{2}}\circ (\hat{e}'_{a; z_{1}-z_{2}}\boxtimes_{P(z_{1})} 
I_{W^{a}})\circ
\mathcal{A}_{P(z_{1}), P(z_{2})}^{P(z_{1}-z_{2}),
P(z_{2})}\circ
(I_{W^{a}}\boxtimes_{P(z_{1})} i'_{a; z_{2}})\circ r_{a; z_{2}}^{-1}.
\end{eqnarray*} 
From this equality and the definitions of $\mu_{a}$ and 
$\mu_{a; z_{1}^{0}, z_{2}^{0}}$,
we obtain
$$\mu_{a}=\mu_{a; z_{1}^{0}, z_{2}^{0}}.$$
By Proposition \ref{main-calc-2}, we obtain (\ref{mu-1}).
\epfv

For the compositions 
\begin{eqnarray*}
\lambda'_{a}&=&r_{a'}\circ (I_{(W^{a})'}\boxtimes \hat{e}'_{a})\circ 
\mathcal{A}^{-1}\circ (i'_{a}\boxtimes I_{W^{a}})\circ l^{-1}_{a'},\\
\mu'_{a}&=&l_{a'}\circ (\hat{e}_{a} \boxtimes I_{(W^{a})'})\circ 
\mathcal{A}\circ (I_{W^{a}} \boxtimes i_{a})\circ r^{-1}_{a'}
\end{eqnarray*}
of the module maps in the 
sequences
\begin{equation}\label{map-3}
\begin{CD}
(W^{a})'@>l^{-1}_{a'}>> V\boxtimes (W^{a})'
@>>> \\ 
@>i'_{a}\boxtimes I_{W^{a}}>>
((W^{a})'\boxtimes W^{a})\boxtimes 
(W^{a})'@>>>  \\
 @>\mathcal{A}^{-1}>> (W^{a})'\boxtimes (W^{a}\boxtimes 
(W^{a})')@>I_{(W^{a})'}\boxtimes \; \hat{e}'_{a}>>
(W^{a})'\boxtimes V @>r_{a'}>> (W^{a})',
\end{CD}
\end{equation}
\begin{equation}\label{map-4}
\begin{CD}
(W^{a})'@>r^{-1}_{a'}>> (W^{a})'\boxtimes V
@>>> \\ 
@>I_{W^{a}} \boxtimes \; i_{a}>>
 (W^{a})'\boxtimes (W^{a}\boxtimes (W^{a})')
@>>> \\ 
 @>\mathcal{A}>> ((W^{a})'\boxtimes W^{a})\boxtimes 
(W^{a})'@>\hat{e}_{a} \boxtimes I_{(W^{a})'}>>
V\boxtimes (W^{a})' @>l_{a'}>> (W^{a})',
\end{CD}
\end{equation}
respectively, we have the following:

\begin{prop}\label{main-cor-3}
For $a\in \mathcal{A}$, we have
\begin{eqnarray}\label{lambda'-mu'}
\lambda'_{a}&=&\mu'_{a}\nn
&=&F_{a}\nn
&=&F(\Y_{ae; 1}^{a} \otimes \Y_{a'a; 1}^{e};
\Y_{ea; 1}^{a}\otimes \Y_{aa'; 1}^{e}).
\end{eqnarray}
\end{prop}

The proof of Proposition \ref{main-cor-3} is completely analogous to the 
proofs of Propositions \ref{main-cor} and \ref{main-cor-2}
and is omitted here.

Since $F_{a}\ne 0$, we let
\begin{eqnarray*}
e_{a}&=&\frac{1}{F_{a}}\hat{e}_{a},\\
e'_{a}&=&\frac{1}{F_{a}}\hat{e}'_{a}
\end{eqnarray*}
for $a\in \A$.

\begin{thm}\label{rigid}
Let $V$ be a simple vertex operator algebra satisfying the 
conditions in Section 1.
Then with the left duals, the right duals and the 
module maps $i_{a}$, $e_{a}$, $i'_{a}$ and $e'_{a}$ above, 
the braided tensor category structure on the 
category of $V$-modules constructed in 
\cite{HL1}--\cite{HL4}, \cite{H1} and \cite{H5}
is rigid.
\end{thm}
\pf 
Since our 
tensor category is semisimple, to prove the rigidity, we need 
only discuss irreducible modules. 
By Propositions \ref{main-cor}, \ref{main-cor-2} and 
\ref{main-cor-3} and the definitions of the maps $e_{a}$ and $e'_{a}$, 
the composition of the maps in the 
sequences
$$\begin{CD}
W^{a}@>l^{-1}_{a}>> V\boxtimes W^{a} @>i_{a}\boxtimes I_{W^{a}}>>
(W^{a}\boxtimes (W^{a})')\boxtimes 
W^{a}@>>>\\
@>\A^{-1}>>W^{a}\boxtimes ((W^{a})'\boxtimes 
W^{a})@>I_{W^{a}}\boxtimes e_{a}>>
W^{a}\boxtimes V@>r_{a}>> W_{a},
\end{CD}$$
$$
\begin{CD}
W^{a}@>r^{-1}_{a}>> W^{a}\boxtimes V
@>I_{W^{a}} \boxtimes \; i'_{a}>>
W^{a}\boxtimes ((W^{a})'\boxtimes W^{a})
@>>>  \\
 @>\mathcal{A}>> (W^{a}\boxtimes (W^{a})')\boxtimes 
W^{a}@>e'_{a} \boxtimes I_{W^{a}}>>
V\boxtimes W^{a} @>l_{a}>> W^{a},
\end{CD}
$$
$$\begin{CD}
(W^{a})'@>l^{-1}_{a'}>> V\boxtimes (W^{a})'
@>>> \\ 
@>i'_{a}\boxtimes I_{W^{a}}>>
((W^{a})'\boxtimes W^{a})\boxtimes 
(W^{a})'@>>>  \\
 @>\mathcal{A}^{-1}>> (W^{a})'\boxtimes (W^{a}\boxtimes 
(W^{a})')@>I_{(W^{a})'}\boxtimes \; e'_{a}>>
(W^{a})'\boxtimes V @>r_{a'}>> (W^{a})',
\end{CD}
$$
$$\begin{CD}
(W^{a})'@>r^{-1}_{a'}>> (W^{a})'\boxtimes V
@>>> \\ 
@>I_{W^{a}} \boxtimes \; i_{a}>>
 (W^{a})'\boxtimes (W^{a}\boxtimes (W^{a})')
@>>> \\ 
 @>\mathcal{A}>> ((W^{a})'\boxtimes W^{a})\boxtimes 
(W^{a})'@>e_{a} \boxtimes I_{(W^{a})'}>>
V\boxtimes (W^{a})' @>l_{a'}>> (W^{a})'
\end{CD}
$$
are the identity maps. So the tensor category is rigid.
\epfv

The calculations above also allow us to calculate the tensor-categorical 
dimensions of $V$-modules:

\begin{thm}
For $a\in \A$, the tensor-categorical dimension of $W^{a}$ 
is $\frac{1}{F_{a}}$.
\end{thm}
\pf
For any $a\in \A$, choose  $w_{a}\in w^{a}$ and
$w'_{a}\in (W^{a})'$ to be lowest weight vectors 
such that $\langle w'_{a}, w_{a}\rangle=1$.
Then by the definitions of $i_{a}$ and $\hat{e}'_{a}$, 
\begin{eqnarray*} 
i_{a}(\mathbf{1})&=&P_{0}(w_{a}\boxtimes w'_{a}),\\
\overline{\hat{e}'_{a}}(w_{a}\boxtimes w'_{a})&=&
\Y_{aa'; 1}^{e}(w_{a}, 1)w'_{a}.
\end{eqnarray*}
Thus we have 
\begin{eqnarray*}
\hat{e}'_{a}(i_{a}(\mathbf{1}))&=&\hat{e}'_{a}(P_{0}(w_{a}\boxtimes w'_{a}))\nn
&=&P_{0}(\overline{\hat{e}'_{a}}(w_{a}\boxtimes w'_{a}))\nn
&=&P_{0}(\Y_{aa'; 1}^{e}(w_{a}, 1)w'_{a})\nn
&=&\res_{x}x^{2h_{a}-1}
\Y_{aa'; 1}^{e}(w_{a}, x)w'_{a}.
\end{eqnarray*}
Since $V$ is irreducible and $V_{(0)}=\C\mathbf{1}$, 
there exists $\nu \in \C$ such that
$$\res_{x}x^{2h_{a}-1}
\Y_{aa'; 1}^{e}(w_{a}, x)w'_{a}=\nu \mathbf{1}.$$
We now calculate $\nu$. 
By the definition of the intertwining operator $\Y_{aa';1}^{e}$
and the assumption that $w_{a}$ and $w'_{a}$ are lowest weight 
vectors and $\langle w'_{a}, w_{a}\rangle=1$,
we have
\begin{eqnarray*}
\lefteqn{\res_{x}x^{2h_{a}-1}(\Y_{aa'; 1}^{e}(w_{a}, x)w'_{a},
\mathbf{1})}\nn
&&=e^{\pi ih_{a}}\res_{x}x^{2h_{a}-1}\langle w'_{a},
\Y_{ae; 1}^{a}(e^{xL(1)}e^{-\pi iL(0)} x^{-2L(0)}w_{a}, x^{-1})\mathbf{1}\rangle\nn
&&=\res_{x}x^{-1}\langle w'_{a}, \Y_{ae; 1}^{a}(w_{a}, x^{-1})\mathbf{1}\rangle\nn
&&=\res_{x}x^{-1}\langle w'_{a}, e^{x^{-1}L(-1)}w_{a}\rangle\nn
&&=\res_{x}x^{-1}\langle  e^{x^{-1}L(1)}w'_{a}, w_{a}\rangle\nn
&&=\langle w'_{a}, w_{a}\rangle\nn
&&=1.
\end{eqnarray*}
Since $(\mathbf{1}, \mathbf{1})=1$, we obtain $\nu=1$. 
Since $V$ is irreducible, the calculation above shows that 
$$\hat{e}'_{a}\circ i_{a}=I_{V}.$$
Thus 
$$e'_{a}\circ i_{a}=\frac{1}{F_{a}}I_{V}.$$
By the definition of the tensor-categorical dimension (see 
\cite{T}), we see that the dimension of 
$W^{a}$ 
is $\frac{1}{F_{a}}$.
\epfv

\begin{rema}\label{v-rigidity}
{\rm In \cite{HL3}, the notion of vertex 
tensor category was introduced and it was proved in 
\cite{HL3} and  \cite{H5} that for a vertex operator algebra
$V$ satisfying the conditions (weaker than) in this paper, 
the category of $V$-modules has a natural structure of 
vertex tensor category. In the proof of the rigidity above, 
we have calculated the compositions of module maps
in the sequences (\ref{map'-z}), (\ref{map'-z-2}) and 
the corresponding sequences (which we have omitted) for the proof of 
(\ref{main-cor-3}). If we replace $\hat{e}_{a; z_{1}, z_{2}}$ and 
$\hat{e}'_{a; z_{1}, z_{2}}$ by 
\begin{eqnarray*}
e_{a; z_{1}, z_{2}}&=&\frac{1}{F_{a}}\hat{e}_{a; z_{1}, z_{2}},\\
e'_{a; z_{1}, z_{2}}&=&\frac{1}{F_{a}}\hat{e}'_{a; z_{1}, z_{2}},
\end{eqnarray*}
these compositions are equal to identity maps.
Actually for any vertex tensor category
whose objects have left and right duals and also have the associated morphisms,
these sequences still make sense. So
we can introduce a notion of rigidity for vertex tensor categories
by requiring the compositions of the morphisms in these sequences 
to be equal to the identity on the object for every object in the 
category. Then our result above shows that 
the vertex tensor category of $V$-modules is rigid. 
We know that a vertex tensor category gives a braided tensor category
\cite{HL3}. Then the proofs of Propositions \ref{main-cor}
and \ref{main-cor-2} and the corresponding parts in the proof
(which we have omitted) 
of Proposition \ref{main-cor-3} shows that 
if a vertex tensor category is rigid, then the corresponding braided 
tensor category is also rigid.}
\end{rema}

\renewcommand{\theequation}{\thesection.\arabic{equation}}
\renewcommand{\thethm}{\thesection.\arabic{thm}}
\setcounter{equation}{0}
\setcounter{thm}{0}

\section{Balancing axioms, the nondegeneracy property and 
the modular tensor category}

In this section we prove that the category of $V$-modules 
for a vertex operator algebra $V$ satisfying the conditions in 
Section 1 has a natural structure of modular tensor category. 
The main work is a proof of  the nondegeneracy property.

We first recall the notion of modular tensor category (see \cite{T}
and \cite{BK} for details).
A {\it ribbon category} is a 
rigid braided tensor category (with tensor product bifunctor $\boxtimes$,
the braiding isomorphism $\mathcal{C}$, the unit object $V$ and 
the right dual functor $*$) together with an isomorphism $\theta_{W}
\in \hom(W, W)$ for 
each object $W$ satisfying the following {\it balancing axioms}: 
(i) $\theta_{W_{1}\boxtimes W_{2}}=\mathcal{C}^{2}
\circ(\theta_{W_{1}}\boxtimes \theta_{W_{2}})$. (ii) $\theta_{V}=I_{V}$. (iii)
$\theta_{W^{*}}=(\theta_{W})^{*}$.
A semisimple ribbon category with 
finitely many inequivalent irreducible objects is a {\it modular 
tensor category} if it has the following {\it nondegeneracy 
property}: Let $\{W_{1}, \dots, W_{m}\}$ be a complete set of representatives 
of equivalence classes of irreducible objects. Then
the $m\times m$ matrix formed by the traces
of the morphism $\mathcal{C}^{2}\in \hom(W_{i}\boxtimes W_{j}, 
W_{i}\boxtimes W_{j})$
in the ribbon category is invertible. 
See, for example, \cite{T} and \cite{BK} for 
details of the theory of modular tensor 
categories. 

Now we consider a vertex operator algebra $V$ satisfying the conditions 
in Section 1.
For a $V$-module $W$, let 
$$\theta_{W}=e^{2\pi iL(0)}.$$

\begin{thm}
The rigid braided tensor category of $V$-modules together with the balancing 
isomorphism or 
the twist $\theta_{W}$ for each $V$-module $W$ is a ribbon category.
\end{thm}
\pf
Let $W_{1}$ and $W_{2}$ be $V$-modules and 
let $\Y$ be the intertwining operator of type ${W_{1}\boxtimes W_{2}\choose 
W_{1}W_{2}}$ such that $w_{1}\boxtimes w_{2}=\Y(w_{1}, 1)w_{2}$ for $w_{1}\in W_{1}$
and $w_{2}\in W_{2}$. From the definition of the braiding isomorphism 
$\mathcal{C}$, we have 
\begin{eqnarray*}
\overline{\mathcal{C}^{2}}(w_{1}\boxtimes w_{2})
&=&\overline{\mathcal{C}^{2}}(
\Y(w_{1}, 1)w_{2})\nn
&=&\Y(w_{1}, x)w_{2}|_{x^{n}=e^{2\pi ni}, \; n\in \R}
\end{eqnarray*}
for $w_{1}\in W_{1}$ and $w_{2}\in W_{2}$. 
Then for $w_{1}\in W_{1}$ and $w_{2}\in W_{2}$, 
\begin{eqnarray*}
\overline{\theta_{W_{1}\boxtimes W_{2}}}(w_{1}\boxtimes w_{2})&=&
e^{2\pi iL(0)}(w_{1}\boxtimes w_{2})\nn
&=&e^{2\pi iL(0)}\Y((w_{1}, 1)w_{2})\nn
&=&\Y(e^{2\pi iL(0)}w_{1}, x)e^{2\pi iL(0)}w_{2}|_{x^{n}=e^{2\pi ni}, \; n\in \R}\nn
&=&\overline{\mathcal{C}^{2}}(
(e^{2\pi iL(0)}w_{1}) \boxtimes (e^{2\pi iL(0)}w_{2}))\nn
&=&\overline{\mathcal{C}^{2}}\circ
\overline{\theta_{W_{1}}\boxtimes \theta_{W_{2}}}(w_{1}\boxtimes w_{2})\nn
&=&\overline{\mathcal{C}^{2}\circ
\theta_{W_{1}}\boxtimes \theta_{W_{2}}}(w_{1}\boxtimes w_{2}).
\end{eqnarray*}
Since the homogeneous components of $w_{1}\boxtimes w_{2}$ for all $w_{1}\in W_{1}$
and $w_{2}\in W_{2}$ span the $V$-module $W_{1}\boxtimes W_{2}$, 
we obtain 
$$\theta_{W_{1}\boxtimes W_{2}}=\mathcal{C}^{2}
\circ(\theta_{W_{1}}\boxtimes \theta_{W_{2}}).$$
Since $V$ is $\Z$-graded, we have $\theta_{V}=I_{V}$. 
For any $V$-module $W$, 
\begin{eqnarray*}
(\theta_{W})^{*}&=&(e^{2\pi iL(0)})^{*}\nn
&=&e^{2\pi iL'(0)}\nn
&=&\theta_{W'}.
\end{eqnarray*}
Thus all the balancing axioms are satisfied.
So the category of $V$-module is a ribbon category.
\epfv

To prove the nondegeneracy property, we calculate the compositions 
of the maps in the sequence
\begin{eqnarray}\label{link}
\begin{CD}
V\\
@Vi'_{a_{2}}VV\\
(W^{a_{2}})'\boxtimes W^{a_{2}}\\
@VI_{(W^{a_{2}})'}\boxtimes r^{-1}_{a_{2}}
VV\\
(W^{a_{2}})'\boxtimes
(W^{a_{2}}\boxtimes V)\\
@VI_{(W^{a_{2}})'}\boxtimes (I_{W^{a_{2}}}\boxtimes\;
i'_{a_{1}; })VV\\
(W^{a_{2}})'\boxtimes
(W^{a_{2}}\boxtimes
((W^{a_{1}})'\boxtimes W^{a_{1}}))\\
@VI_{(W^{a_{2}})'}\boxtimes
\mathcal{A} VV\\
(W^{a_{2}})'\boxtimes
((W^{a_{2}}\boxtimes
(W^{a_{1}})')\boxtimes W^{a_{1}})\\
@VI_{(W^{a_{2}})'}\boxtimes
(\mathcal{C}^{2}\boxtimes
I_{W^{a_{2}}})VV\\ 
(W^{a_{2}})'\boxtimes
((W^{a_{2}}\boxtimes
(W^{a_{1}})')\boxtimes W^{a_{1}})\\
@VI_{(W^{a_{2}})'}\boxtimes
\mathcal{A}VV\\
(W^{a_{2}})'\boxtimes
(W^{a_{2}}\boxtimes
((W^{a_{1}})'\boxtimes W^{a_{1}}))\\
@VI_{(W^{a_{2}})'}\boxtimes
(I_{W^{a_{2}}}\boxtimes
\hat{e}_{a_{1}})VV\\
(W^{a_{2}})'\boxtimes
(W^{a_{2}}\boxtimes
V)\\
@VI_{(W^{a_{2}})'}\boxtimes
r_{a_{2}}VV\\
(W^{a_{2}})'\boxtimes W^{a_{2}}\\
@V\hat{e}_{a_{2}}VV\\
V
\end{CD}&\nn
&&
\end{eqnarray}
This composition will give us the trace of $\mathcal{C}^{2}\in \hom((W^{a_{1}})', 
W^{a_{2}})$.
Similarly to the calculations of the maps in the proof of the rigidity
in the preceding section, we first calculate the composition 
of certain other module maps.

For $z_{1}, z_{2}, z_{3}\in \C^{\times}$ satisfying $|z_{1}|>|z_{2}|
>|z_{3}|>|z_{2}-z_{3}|>0$,
consider the module map given by the sequence
\begin{eqnarray}\label{v-link}
&\begin{CD}
V\\
@Vi'_{a_{2}; z_{1}}VV\\
(W^{a_{2}})'\boxtimes_{P(z_{1})} W^{a_{2}}\\
@VI_{(W^{a_{2}})'}\boxtimes_{P(z_{1})}r^{-1}_{a_{2}, z_{2}}
VV\\
(W^{a_{2}})'\boxtimes_{P(z_{1})} 
(W^{a_{2}}\boxtimes_{P(z_{2})} V)\\
@VI_{(W^{a_{2}})'}\boxtimes_{P(z_{1})}(I_{W^{a_{2}}}\boxtimes_{P(z_{2})}
i'_{a_{1}; z_{3}})VV\\
(W^{a_{2}})'\boxtimes_{P(z_{1})} 
(W^{a_{2}}\boxtimes_{P(z_{2})} 
((W^{a_{1}})'\boxtimes_{P(z_{3})} W^{a_{1}}))\\
@VI_{(W^{a_{2}})'}\boxtimes_{P(z_{1})} 
\mathcal{A}_{P(z_{2}), P(z_{2})}
^{P(z_{2}-z_{3}), P(z_{3})}VV\\
(W^{a_{2}})'\boxtimes_{P(z_{1})} 
((W^{a_{2}}\boxtimes_{P(z_{2}-z_{3})} 
(W^{a_{1}})')\boxtimes_{P(z_{3})} W^{a_{1}})\\
@VI_{(W^{a_{2}})'}\boxtimes_{P(z_{1})} 
((\mathcal{C}_{P(z_{2}-z_{3})})^{2}\boxtimes_{P(z_{3})} 
I_{W^{a_{2}}})VV\\ 
(W^{a_{2}})'\boxtimes_{P(z_{1})} 
((W^{a_{2}}\boxtimes_{P(z_{2}-z_{3})} 
(W^{a_{1}})')\boxtimes_{P(z_{3})} W^{a_{1}})\\
@VI_{(W^{a_{2}})'}\boxtimes_{P(z_{1})} 
(\mathcal{A}_{P(z_{2}), P(z_{2})}
^{P(z_{2}-z_{3}), P(z_{3})})^{-1}VV\\
(W^{a_{2}})'\boxtimes_{P(z_{1})} 
(W^{a_{2}}\boxtimes_{P(z_{2})} 
((W^{a_{1}})'\boxtimes_{P(z_{3})} W^{a_{1}}))\\
@VI_{(W^{a_{2}})'}\boxtimes_{P(z_{1})}(I_{W^{a_{2}}}\boxtimes_{P(z_{2})}
\hat{e}_{a_{1}; z_{3}})VV\\(W^{a_{2}})'\boxtimes_{P(z_{1})} 
(W^{a_{2}}\boxtimes_{P(z_{2})}
V)\\
@VI_{(W^{a_{2}})'}\boxtimes_{P(z_{1})}
r_{a_{2}; z_{2}}VV\\(W^{a_{2}})'\boxtimes_{P(z_{1})} W^{a_{2}}\\
@V\hat{e}_{a_{2}; z_{1}}VV\\V
\end{CD}&\nn
&&
\end{eqnarray}

\begin{prop}\label{v-link-calc-1}
The composition of the module maps in the sequence (\ref{v-link})
is equal to $(B^{(-1)})^{2}_{a_{2}, a_{1}}I_{V}.$
\end{prop}
\pf
Since $V$ is an irreducible $V$-module, we know that 
the composition must be proportional to $I_{V}$. So we need
only show that the proportional constant is 
$(B^{(-1)})^{2}_{a_{2}, a_{1}}$. To show this, we need 
only calculate the composition applied to any 
nonzero element of $V$ or the natural extension of the composition 
applied to any nonzero element of $\overline{V}$. 

By definition, we have
$$w'\boxtimes_{P(z)}w=\overline{(f'_{a; z})^{-1}}\left(
\sum_{b\in \A}\sum_{k=1}^{N_{a'a}^{b}}
\Y_{a'a; k}^{b}(w', z)w\right)$$
for $a\in \A$, $z\in \C$, $w\in W^{a}$ and
$w'\in (W^{a})'$. Then for any module maps
$\alpha: (W^{a})'\to (W^{a})'$ and 
$\beta: W^{a}\to W^{a}$, 
\begin{eqnarray}\label{f'-commu}
\lefteqn{(\overline{\alpha\boxtimes_{P(z)}\beta})
\left(\overline{(f'_{a; z})^{-1}}
\left(\sum_{b\in \A}\sum_{k=1}^{N_{a'a}^{b}}
\Y_{a'a; k}^{b}(w', z)w\right)\right)}\nn
&&=(\overline{\alpha\boxtimes_{P(z)}\beta})(w'\boxtimes_{P(z)}w)\nn
&&=\alpha(w')\boxtimes_{P(z)}\beta(w)\nn
&&=\overline{(f'_{a; z})^{-1}}\left(
\sum_{b\in \A}\sum_{k=1}^{N_{a'a}^{b}}
\Y_{a'a; k}^{b}(\alpha(w'), z)\beta(w)\right).
\end{eqnarray}
Since $f'_{a; z}$ is an equivalence of $V$-modules and 
the homogeneous components of elements of the form
$w'\boxtimes_{P(z)}w$ for $w\in W^{a}$ and
$w'\in (W^{a})'$
span $(W^{a})'\boxtimes_{P(z)}W^{a}$,  
the homogeneous components of elements of the form
\begin{equation}\label{span}
\sum_{b\in \A}\sum_{k=1}^{N_{a'a}^{b}}
\Y_{a'a; k}^{b}(w', z)w
\end{equation}
for $w\in W^{a}$ and
$w'\in (W^{a})'$ span $\coprod_{b\in \A}N_{a'a}^{b}W^{b}$.
In particular, by the definition of $g'_{a; z}$, 
for $w_{a}\in W^{a}$ and
$w'_{a}\in (W^{a})'$, $\overline{g'_{a; z}}
(\Y_{a'a; 1}^{e}(w'_{a}, z)w_{a})$
can be obtained as a sum of homogeneous components of 
elements of the form (\ref{span}) and 
$\overline{g'_{a; z}}(\Y_{a'a; 1}^{e}(\alpha(w'_{a}), z)\beta(w_{a}))$
is the sum of the same homogeneous components of the same
elements of the form (\ref{span}) with $w'$ and $w$ replaced 
by $\alpha(w')$ and $\beta(w)$, respectively. Thus from 
(\ref{f'-commu}) we obtain
\begin{eqnarray}\label{f'-g'-commu}
\lefteqn{(\overline{\alpha\boxtimes_{P(z)}\beta})
(\overline{(f'_{a; z})^{-1}}(
\overline{g'_{a; z}}(\Y_{a'a; 1}^{e}(w'_{a}, z)w_{a})))}\nn
&&=\overline{(f'_{a; z})^{-1}}(
\overline{g'_{a; z}}(\Y_{a'a; 1}^{e}(\alpha(w'_{a}), z)\beta(w_{a}))).
\end{eqnarray}
From the definition of $i'_{a; z}$ and (\ref{f'-g'-commu}), 
we obtain
\begin{equation}\label{i'-commu}
(\overline{\alpha\boxtimes_{P(z)} \beta)\circ i'_{a; z}})
(\Y_{a'a; 1}^{e}(w'_{a}, z)w_{a})
=\overline{i'_{a; z}} (\Y_{a'a; 1}^{e}(\alpha(w'_{a}), z)\beta(w_{a}))
\end{equation}
for module maps $\alpha: (W^{a})'\to (W^{a})'$ and 
$\beta: W^{a}\to W^{a}$,  $a\in \A$, $z\in \C$, $w_{a}\in W^{a}$ and
$w'_{a}\in (W^{a})'$.
Using (\ref{i'-commu}),  
we obtain
\begin{eqnarray}\label{modular-calc1}
\lefteqn{(\overline{I_{(W^{a_{2}})'}\boxtimes_{P(z_{1})}
(I_{W^{a_{2}}}\boxtimes_{P(z_{2})}
i'_{a_{1}; z_{3}})})
\biggl((\overline{I_{(W^{a_{2}})'}\boxtimes_{P(z_{1})}
r^{-1}_{a_{2}, z_{2}}})}\nn
&&\biggl(\overline{i'_{a_{2}; z_{1}}}
\biggl(\Y_{a'_{2}a_{2}; 1}^{e}(w'_{a_{2}}, 
z_{1})\Y_{a_{2}e; 1}^{a_{2}}
(w_{a_{2}}, z_{2})\Y_{a_{1}'a_{1}; 1}^{e}(w_{a'_{1}}, z_{3})
w_{a_{1}}\biggr)\biggr)\biggr)\nn
&&=\overline{i'_{a_{2}; z_{1}}}
\biggl(\Y_{a'_{2}a_{2}; 1}^{e}(w'_{a_{2}}, z_{1})\nn
&&\quad\quad\quad\quad\quad
\biggl(w_{a_{2}}\boxtimes_{P(z_{2})}
\biggl(\overline{i'_{a_{1}; z_{3}}}
\biggl(\Y_{a_{1}'a_{1}; 1}^{e}(w_{a'_{1}}, z_{3})
w_{a_{1}}\biggr)\biggr)\biggr)\biggr)
\end{eqnarray}
for $w_{a_{1}}\in W^{a_{1}}$, $w_{a_{2}}\in W^{a_{2}}$,
$w'_{a_{1}}\in (W^{a_{1}})'$ and $w'_{a_{2}}\in (W^{a_{2}})'$.

We also have the following commutative diagram
\begin{eqnarray}\label{braiding-diag}
\begin{CD}
W^{a_{2}}\boxtimes_{P(z_{2})}
((W^{a_{1}})'\boxtimes_{P(z_{3})}
W^{a_{1}})@>r_{a_{2}; z_{2}}\circ (I_{W^{a_{2}}}\boxtimes_{P(z_{2})} 
\hat{e}_{a'_{1}; z_{3}})>>W^{a_{2}}\\
@V\mathcal{M}VV @VVm_{a_{2}, a_{1}}V\\
W^{a_{2}}\boxtimes_{P(z_{2})}
((W^{a_{1}})'\boxtimes_{P(z_{3})}
W^{a_{1}})@>r_{a_{2}; z_{2}}\circ (I_{W^{a_{2}}}\boxtimes_{P(z_{2})} 
\hat{e}_{a'_{1}; P(z_{3})})>>W^{a_{2}}
\end{CD}\nn
&&
\end{eqnarray}
where 
$$\mathcal{M}=
\left(\mathcal{A}_{P(z_{2}), P(z_{2})}^{P(z_{2}-z_{3}), 
P(z_{3})}\right)^{-1}
\circ ((\mathcal{C}_{P(z_{2}-z_{3})})^{2}
\boxtimes_{P(z_{3})} 
I_{W^{a_{2}}})\circ 
\mathcal{A}_{P(z_{2}), P(z_{2})}^{P(z_{2}-z_{3}), 
P(z_{3})}$$
and the module map $m_{a_{2}, a_{1}}$ is defined as follows:
Consider the element 
\begin{equation}\label{element}
\Y_{a_{2}e; 1}^{a_{2}}
(w_{a_{2}}, z_{2})\Y_{a_{1}'a_{1}; 1}^{e}(w'_{a_{1}}, z_{3})
w_{a_{1}}
\end{equation}
of $\overline{W^{a_{2}}}$ and the path $\gamma$ in 
$M^{2}=\{(z_{2}, z_{3})\in \C^{2}\;|\;  z_{2}, z_{3}\ne 0, 
z_{2}\ne z_{3}\}$ given by 
$$\gamma(t)=(z_{3}+e^{-2\pi it}(z_{2}-z_{3}), z_{3}).$$
Starting with the element (\ref{element}) above, the analytic extension
along the path $\gamma$ gives a value at the endpoint of  $\gamma$.
This is again an element of $\overline{W^{a_{2}}}$ which we denote 
$$\overline{m_{a_{2}, a_{1}}}(\Y_{a_{2}e; 1}^{a_{2}}
(w_{a_{2}}, z_{2})\Y_{a_{1}'a_{1}; 1}^{e}(w'_{a_{1}}, z_{3})
w_{a_{1}}).$$ 
The correspondence 
\begin{eqnarray*}
\lefteqn{\Y_{a_{2}e; 1}^{a_{2}}
(w_{a_{2}}, z_{2})\Y_{a_{1}'a_{1}; 1}^{e}(w'_{a_{1}}, z_{3})
w_{a_{1}}}\nn
&&
\mapsto \overline{m_{a_{2}, a_{1}}}(\Y_{a_{2}e; 1}^{a_{2}}
(w_{a_{2}}, z_{2})\Y_{a_{1}'a_{1}; 1}^{e}(w'_{a_{1}}, z_{3})
w_{a_{1}})
\end{eqnarray*}
determines a unique module map $m_{a_{2}, a_{1}}$ from 
$W^{a_{2}}$ to itself. 

Note that by definition, the module map $\mathcal{M}$
can also be directly obtained in the same way as that for $m_{a_{2}, a_{1}}$: 
Starting with the element 
$$w_{a_{2}}\boxtimes_{P(z_{2})} 
(w_{a_{1}}' \boxtimes_{P(z_{3})}w_{a_{1}}),$$
the analytic extension 
along the path $\gamma$ gives a value at the endpoint of $\gamma$.
This value is 
$$\overline{\mathcal{M}}(w_{a_{2}}\boxtimes_{P(z_{2})} 
(w_{a_{1}}' \boxtimes_{P(z_{3})}w_{a_{1}}))$$
and the correspondence
$$w_{a_{2}}\boxtimes_{P(z_{2})} 
(w_{a_{1}}' \boxtimes_{P(z_{3})}w_{a_{1}})\mapsto 
\overline{\mathcal{M}}(w_{a_{2}}\boxtimes_{P(z_{2})} 
(w_{a_{1}}' \boxtimes_{P(z_{3})}w_{a_{1}}))$$
determines the map $\mathcal{M}$.
From this construction of $\mathcal{M}$, we see that the 
diagram (\ref{braiding-diag}) is commutative since
analytic extensions certainly commute with the module map 
$$r_{a_{2}; z_{2}}\circ (I_{W^{a_{2}}}\boxtimes_{P(z_{2})} 
e_{a'_{1}; z_{3}}).$$

Now this commutativity of (\ref{braiding-diag}) gives
\begin{eqnarray}\label{braiding-calc0}
\lefteqn{r_{a_{2}; z_{2}}\circ 
(I_{W^{a_{2}}}\boxtimes_{P(z_{2})} 
e_{a'_{1}; z_{3}})\circ 
\left(\mathcal{A}_{P(z_{2}, 
P(z_{2})}^{P(z_{2}-z_{3}), 
P(z_{3})}\right)^{-1}}\nn
&&\quad\quad\quad\quad
\circ ((\mathcal{C}_{P(z_{2}-z_{3})})^{2}
\boxtimes_{P(z_{3})} 
I_{W^{a_{2}}})\circ 
\mathcal{A}_{P(z_{2}, P(z_{2})}^{P(z_{2}-z_{3}), 
P(z_{3})}\nn
&&=m_{a_{2}, a_{1}}\circ 
(r_{a_{2}; z_{2}}\circ (I_{W^{a_{2}}}\boxtimes_{P(z_{2})} 
e_{a'_{1}; z_{3}})).
\end{eqnarray}
On the other hand,
\begin{eqnarray}\label{braiding-calc1}
\overline{m_{a_{2}, a_{1}}\circ 
(r_{a_{2}; z_{2}}\circ (I_{W^{a_{2}}}\boxtimes_{P(z_{2})} 
e_{a'_{1}; z_{3}}))}(w_{a_{2}}\boxtimes_{P(z_{2})} 
(w_{a_{1}}' \boxtimes_{P(z_{3})}w_{a_{1}}))\nn
=\overline{m_{a_{2}, a_{1}}}(\Y_{a_{2}e; 1}^{a_{2}}
(w_{a_{2}}, z_{2})\Y_{a_{1}'a_{1}; 1}^{e}(w'_{a_{1}}, z_{3})
w_{a_{1}}).\nn
\end{eqnarray}
But by definition, for any $w_{a_{2}}'\in (W^{a_{2}})'$, we have
\begin{eqnarray*}
\lefteqn{\langle w_{a_{2}}', 
\overline{m_{a_{2}, a_{1}}}(\Y_{a_{2}e; 1}^{a_{2}}
(w_{a_{2}}, z_{2})\Y_{a_{1}'a_{1}; 1}^{e}(w'_{a_{1}}, z_{3})
w_{a_{1}})\rangle}\nn
&&=(B^{(-1)})^{2}(\langle w_{a_{2}}', 
\Y_{a_{2}e; 1}^{a_{2}}
(w_{a_{2}}, z_{2})\Y_{a_{1}'a_{1}; 1}^{e}(w'_{a_{1}}, z_{3})
w_{a_{1}}\rangle).
\end{eqnarray*}
Thus by 
the commutativity of intertwining operators, we have
\begin{eqnarray}\label{braiding-calc2}
\lefteqn{\overline{m_{a_{2}, a_{1}}}(\Y_{a_{2}e; 1}^{a_{2}}
(w_{a_{2}}, z_{2})\Y_{a_{1}'a_{1}; 1}^{e}(w'_{a_{1}}, z_{3})
w_{a_{1}})}\nn
&&=(B^{(-1)})^{2}(\Y_{a_{2}e; 1}^{a_{2}}\otimes 
\Y_{a_{1}'a_{1}; 1}^{e}; \Y_{a_{2}e; 1}^{a_{2}}\otimes
\Y_{a_{1}'a_{1}; 1}^{e})\cdot \nn
&&\quad\quad\quad\quad\quad\quad\quad\quad\quad \cdot\Y_{a_{2}e; 1}^{a_{2}}
(w_{a_{2}}, z_{2})\Y_{a_{1}'a_{1}; 1}^{e}(w'_{a_{1}}, z_{3})
w_{a_{1}}.
\end{eqnarray}
Combining (\ref{braiding-calc0})--(\ref{braiding-calc2}),
we obtain
\begin{eqnarray}\label{braiding-calc3}
\lefteqn{\overline{(r_{a_{2}; z_{2}}}
\biggl(\overline{(I_{W^{a_{2}}}\boxtimes_{P(z_{2})} 
\hat{e}_{a'_{1}; z_{3}}))}}\nn
&&\quad\quad\quad\quad\biggl(\overline{
\left(\mathcal{A}_{P(z_{2}), 
P(z_{2})}^{P(z_{2}-z_{3}), 
P(z_{3})}\right)^{-1}}
\biggl(\overline{((\mathcal{C}_{P(z_{2}-z_{3})})^{2}
\boxtimes_{P(z_{3})} 
I_{W^{a_{2}}})}\nn
&&\quad\quad\quad\quad\quad\quad\quad\quad\biggl(\overline{
\mathcal{A}_{P(z_{2}), P(z_{2})}^{P(z_{2}-z_{3}), 
P(z_{3})}}
(w_{a_{2}}\boxtimes_{P(z_{2})} 
(w'_{a_{1}} \boxtimes_{P(z_{3})}w_{a_{1}}))\biggr)\biggr)\biggr)
\biggr)\nn
&&=(B^{(-1)})^{2}(\Y_{a_{2}e; 1}^{a_{2}}\otimes 
\Y_{a_{1}'a_{1}; 1}^{e}; \Y_{a_{2}e; 1}^{a_{2}}\otimes
\Y_{a_{1}'a_{1}; 1}^{e})\cdot \nn
&&\quad\quad\quad\quad\quad\quad\quad\quad\quad \cdot
\Y_{a_{2}e; 1}^{a_{2}}
(w_{a_{2}}, z_{2})\Y_{a_{1}'a_{1}; 1}^{e}(w'_{a_{1}}, z_{3})
w_{a_{1}}.
\end{eqnarray}

Since $\overline{i'_{a_{1}; z_{3}}}
(\Y_{a_{1}'a_{1}; 1}^{e}(w_{a'_{1}}, z_{3})
w_{a_{1}})$ is in $\overline{(W^{a})'\boxtimes_{P(z_{3})}
W^{a}}$, it can be written as a sum of of homogeneous components of 
elements of the form $w'\boxtimes_{P(z_{3})}w$ for $w'\in (W^{a})'$ 
and $w\in W^{a}$. Since $f'_{a; z_{3}}$ is an equivalence of 
$V$-modules, $g'_{a; z_{3}}(\Y_{a_{1}'a_{1}; 1}^{e}(w_{a'_{1}}, z_{3})
w_{a_{1}})$ is the same sum of the same homogeneous 
components of the elements of the form 
$$\sum_{b\in \A}
\sum_{k=1}^{N_{a_{1}'a_{1}}^{b}}\Y_{a_{1}'a_{1}; k}^{b}(w_{a'_{1}}, z_{3})
w_{a_{1}}.$$
Then we see that the same sum of the same homogeneous 
components of the elements of the form 
$Y_{a_{1}'a_{1}; 1}^{e}(w_{a'_{1}}, z_{3})
w_{a_{1}}$ must be equal to itself. Using these facts and taking the 
coefficients of $z_{3}$ to all powers in both sides of 
(\ref{braiding-calc3}) and then taking the sum above,  we obtain
\begin{eqnarray}\label{braiding-calc4}
\lefteqn{\overline{r_{a_{2}; z_{2}}}
\biggl((\overline{I_{W^{a_{2}}}\boxtimes_{P(z_{2})} 
\hat{e}_{a'_{1}; z_{3}}})}\nn
&&\quad\quad\quad\biggl(
\overline{\left(\mathcal{A}_{P(z_{2}), 
P(z_{2})}^{P(z_{2}-z_{3}), 
P(z_{3})}\right)^{-1}}
\biggl((\overline{(\mathcal{C}_{P(z_{2}-z_{3})})^{2}
\boxtimes_{P(z_{3})} 
I_{W^{a_{2}}}})\nn
&&\quad\quad\quad\quad\quad\quad\biggl(
\overline{\mathcal{A}_{P(z_{2}), P(z_{2})}^{P(z_{2}-z_{3}), 
P(z_{3})}}
(w_{a_{2}}\boxtimes_{P(z_{2})} \overline{i'_{a_{1}; z_{3}}}
(\Y_{a_{1}'a_{1}; 1}^{e}(w_{a'_{1}}, z_{3})
w_{a_{1}}))\biggr)\biggr)\biggr)
\biggr)\nn
&&=(B^{(-1)})^{2}(\Y_{a_{2}e; 1}^{a_{2}}\otimes 
\Y_{a_{1}'a_{1}; 1}^{e}; \Y_{a_{2}e; 1}^{a_{2}}\otimes
\Y_{a_{1}'a_{1}; 1}^{e})\cdot \nn
&&\quad\quad\quad\quad\quad\quad\quad\quad\quad \cdot
\Y_{a_{2}e; 1}^{a_{2}}
(w_{a_{2}}, z_{2})\Y_{a_{1}'a_{1}; 1}^{e}(w'_{a_{1}}, z_{3})
w_{a_{1}}.
\end{eqnarray}

From (\ref{braiding-calc4}),
we obtain
\begin{eqnarray}\label{modular-calc1.5}
\lefteqn{\biggl(\overline{\hat{e}_{a_{2}'; z_{1}}}\circ 
(\overline{I_{(W^{a_{2}})'}\boxtimes_{P(z_{1})}
r_{a_{2}; z_{2}}})
\circ (\overline{I_{(W^{a_{2}})'}\boxtimes_{P(z_{1})}
(I_{W^{a_{2}}}\boxtimes_{P(z_{2})} 
\hat{e}_{a'_{1}; z_{3}}))}}\nn
&&\circ \left(\overline{I_{(W^{a_{2}})'}\boxtimes_{P(z_{1})}
\left(\mathcal{A}_{P(z_{2}), 
P(z_{2})}^{P(z_{2}-z_{3}), 
P(z_{3})}\right)^{-1}}\right)\nn
&&\circ (\overline{I_{(W^{a_{2}})'}\boxtimes_{P(z_{1})}
((\mathcal{C}_{P(z_{2}-z_{3})})^{2}
\boxtimes_{P(z_{3})} 
I_{W^{a_{2}}})})
\circ \left(\overline{I_{(W^{a_{2}})'}\boxtimes_{P(z_{1})}
\mathcal{A}_{P(z_{2}), P(z_{2})}^{P(z_{2}-z_{3}), 
P(z_{3})}}\right)\nn
&&\quad\biggl(w'_{a_{2}}\boxtimes_{P(z_{1})}\biggl(
w_{a_{2}}\boxtimes_{P(z_{2})}\overline{i'_{a_{1}; z_{3}}}
\biggl(\Y_{a_{1}'a_{1}; 1}^{e}(w'_{a_{1}}, z_{3})
w_{a_{1}}\biggr)\biggr)\biggr)\nn
&&=(B^{(-1)})^{2}(\Y_{a_{2}e; 1}^{a_{2}}\otimes 
\Y_{a_{1}'a_{1}; 1}^{e}; \Y_{a_{2}e; 1}^{a_{2}}\otimes
\Y_{a_{1}'a_{1}; 1}^{e})\cdot \nn
&&\quad\quad\quad\quad\quad\quad \cdot
\biggl(\Y_{a'_{2}a_{2}; 1}^{e}(w'_{a_{2}}, 
z_{1})\Y_{a_{2}e; 1}^{a_{2}}
(w_{a_{2}}, z_{2})\Y_{a_{1}'a_{1}; 1}^{e}(w'_{a_{1}}, 
z_{3})
w_{a_{1}}\biggr)
\end{eqnarray}
Noticing that (\ref{modular-calc1.5}) holds for 
all $w_{a_{1}}\in W^{a_{1}}$, $w_{a_{2}}\in W^{a_{2}}$,
$w'_{a_{1}}\in (W^{a_{1}})'$ and $w'_{a_{2}}\in (W^{a_{2}})'$,
and using the same arguments as we have used to obtain 
(\ref{braiding-calc4}) from (\ref{braiding-calc3}), 
from (\ref{modular-calc1.5}), we obtain
\begin{eqnarray}\label{modular-calc1.6}
\lefteqn{\biggl(\overline{\hat{e}_{a_{2}'; z_{1}}}\circ 
(\overline{I_{(W^{a_{2}})'}\boxtimes_{P(z_{1})}
r_{a_{2}; z_{2}}})
\circ (\overline{I_{(W^{a_{2}})'}\boxtimes_{P(z_{1})}
(I_{W^{a_{2}}}\boxtimes_{P(z_{2})} 
\hat{e}_{a'_{1}; z_{3}}))}}\nn
&&\circ \left(\overline{I_{(W^{a_{2}})'}\boxtimes_{P(z_{1})}
\left(\mathcal{A}_{P(z_{2}), 
P(z_{2})}^{P(z_{2}-z_{3}), 
P(z_{3})}\right)^{-1}}\right)\nn
&&\circ (\overline{I_{(W^{a_{2}})'}\boxtimes_{P(z_{1})}
((\mathcal{C}_{P(z_{2}-z_{3})})^{2}
\boxtimes_{P(z_{3})} 
I_{W^{a_{2}}})})
\circ \left(\overline{I_{(W^{a_{2}})'}\boxtimes_{P(z_{1})}
\mathcal{A}_{P(z_{2}), P(z_{2})}^{P(z_{2}-z_{3}), 
P(z_{3})}}\right)\nn
&&\quad\biggl(\overline{i_{a'_{2}; z_{1}}}\biggl(
\Y_{a_{2}'a_{2}; 1}^{e}(w'_{a_{2}}, z_{1})\biggl(
w_{a_{2}}\boxtimes_{P(z_{2})}
\overline{i'_{a_{1}; z_{3}}}
\biggl(\Y_{a_{1}'a_{1}; 1}^{e}(w'_{a_{1}}, z_{3})
w_{a_{1}}\biggr)\biggr)\biggr)\biggr)\nn
&&=(B^{(-1)})^{2}(\Y_{a_{2}e; 1}^{a_{2}}\otimes 
\Y_{a_{1}'a_{1}; 1}^{e}; \Y_{a_{2}e; 1}^{a_{2}}\otimes
\Y_{a_{1}'a_{1}; 1}^{e})\cdot \nn
&&\quad\quad\quad\quad\quad\quad \cdot
\biggl(\Y_{a'_{2}a_{2}; 1}^{e}(w'_{a_{2}}, 
z_{1})\Y_{a_{2}e; 1}^{a_{2}}
(w_{a_{2}}, z_{2})\Y_{a_{1}'a_{1}; 1}^{e}(w'_{a_{1}}, 
z_{3})
w_{a_{1}}\biggr).
\end{eqnarray}
Combining (\ref{modular-calc1}) and (\ref{modular-calc1.6}), 
we obtain
\begin{eqnarray}\label{modular-calc1.7}
\lefteqn{\biggl(\overline{\hat{e}_{a_{2}'; z_{1}}}\circ 
(\overline{I_{(W^{a_{2}})'}\boxtimes_{P(z_{1})}
r_{a_{2}; z_{2}}})
\circ (\overline{I_{(W^{a_{2}})'}\boxtimes_{P(z_{1})}
(I_{W^{a_{2}}}\boxtimes_{P(z_{2})} 
\hat{e}_{a'_{1}; z_{3}}))}}\nn
&&\circ \left(\overline{I_{(W^{a_{2}})'}\boxtimes_{P(z_{1})}
\left(\mathcal{A}_{P(z_{2}), 
P(z_{2})}^{P(z_{2}-z_{3}), 
P(z_{3})}\right)^{-1}}\right)\nn
&&\circ (\overline{I_{(W^{a_{2}})'}\boxtimes_{P(z_{1})}
((\mathcal{C}_{P(z_{2}-z_{3})})^{2}
\boxtimes_{P(z_{3})} 
I_{W^{a_{2}}})})
\circ \left(\overline{I_{(W^{a_{2}})'}\boxtimes_{P(z_{1})}
\mathcal{A}_{P(z_{2}), P(z_{2})}^{P(z_{2}-z_{3}), 
P(z_{3})}}\right)\nn
&&\circ
(\overline{I_{(W^{a_{2}})'}\boxtimes_{P(z_{1})}
(I_{W^{a_{2}}}\boxtimes_{P(z_{2})}
i_{a'_{1}; z_{3}})})
\circ (\overline{I_{(W^{a_{2}})'}\boxtimes_{P(z_{1})}
r^{-1}_{a_{2}, z_{2}}})\circ \overline{i_{a'_{2}; z_{1}}}\biggr)
\nn
&&\quad\biggl(\Y_{a_{2}'a_{2}; 1}^{e}(w'_{a_{2}}, 
z_{1})\Y_{a_{2}e; 1}^{a_{2}}
(w_{a_{2}}, z_{2})
\biggl(\Y_{a_{1}'a_{1}; 1}^{e}(w'_{a_{1}}, z_{3})
w_{a_{1}}\biggr)\biggr)\nn
&&=(B^{(-1)})^{2}(\Y_{a_{2}e; 1}^{a_{2}}\otimes 
\Y_{a_{1}'a_{1}; 1}^{e}; \Y_{a_{2}e; 1}^{a_{2}}\otimes
\Y_{a_{1}'a_{1}; 1}^{e})\cdot \nn
&&\quad\quad\quad\quad\quad\quad \cdot
\biggl(\Y_{a'_{2}a_{2}; 1}^{e}(w'_{a_{2}}, 
z_{1})\Y_{a_{2}e; 1}^{a_{2}}
(w_{a_{2}}, z_{2})\Y_{a_{1}'a_{1}; 1}^{e}(w'_{a_{1}}, 
z_{3})
w_{a_{1}}\biggr)
\end{eqnarray}
Since we can always choose 
$w_{a_{1}}\in W^{a_{1}}$, $w_{a_{2}}\in W^{a_{2}}$,
$w'_{a_{1}}\in (W^{a_{1}})'$ and $w'_{a_{2}}\in (W^{a_{2}})'$
such that 
$$\Y_{a'_{2}a_{2}; 1}^{e}(w'_{a_{2}}, 
z_{1})\Y_{a_{2}e; 1}^{a_{2}}
(w_{a_{2}}, z_{2})\Y_{a_{1}'a_{1}; 1}^{e}(w'_{a_{1}}, 
z_{3})
w_{a_{1}}\ne 0,$$ 
(\ref{modular-calc1.7}) 
gives
\begin{eqnarray}\label{modular-calc2}
\lefteqn{\biggl(\hat{e}_{a_{2}'; z_{1}}\circ 
(I_{(W^{a_{2}})'}\boxtimes_{P(z_{1})}
r_{a_{2}; z_{2}})
\circ (I_{(W^{a_{2}})'}\boxtimes_{P(z_{1})}
(I_{W^{a_{2}}}\boxtimes_{P(z_{2})} 
\hat{e}_{a'_{1}; z_{3}}))}\nn
&&\circ \left(I_{(W^{a_{2}})'}\boxtimes_{P(z_{1})}
\left(\mathcal{A}_{P(z_{2}), 
P(z_{2})}^{P(z_{2}-z_{3}), 
P(z_{3})}\right)^{-1}\right)\nn
&&\circ (I_{(W^{a_{2}})'}\boxtimes_{P(z_{1})}
((\mathcal{C}_{P(z_{2}-z_{3})})^{2}
\boxtimes_{P(z_{3})} 
I_{W^{a_{2}}}))
\circ \left(I_{(W^{a_{2}})'}\boxtimes_{P(z_{1})}
\mathcal{A}_{P(z_{2}), P(z_{2})}^{P(z_{2}-z_{3}), 
P(z_{3})}\right)\nn
&&\circ
(I_{(W^{a_{2}})'}\boxtimes_{P(z_{1})}
(I_{W^{a_{2}}}\boxtimes_{P(z_{2})}
i_{a'_{1}; z_{3}}))
\circ (I_{(W^{a_{2}})'}\boxtimes_{P(z_{1})}
r^{-1}_{a_{2}, z_{2}})\circ i_{a'_{2}; z_{1}}\biggr)
\nn
&&=(B^{(-1)})^{2}(\Y_{a_{2}e; 1}^{a_{2}}\otimes 
\Y_{a_{1}'a_{1}; 1}^{e}; \Y_{a_{2}e; 1}^{a_{2}}\otimes
\Y_{a_{1}'a_{1}; 1}^{e})I_{V}\nn
&&=(B^{(-1)})_{a_{2}, a_{1}}^{2}I_{V}.
\end{eqnarray}
\epfv

\begin{prop}\label{v-link-calc-2}
The composition of the module maps in the sequence (\ref{link})
is equal to $(B^{(-1)})^{2}_{a_{2}, a_{1}}I_{V}.$
\end{prop}
\pf
Choose $z_{1}^{0}, z_{2}^{0}, z_{3}^{0}\in (0, \infty$ satisfying 
$z_{1}^{0}>z_{2}^{0}>z_{3}^{0}>z_{2}^{0}-z_{3}^{0}>0$. Let 
$\gamma_{1}, \gamma_{2}, \gamma_{3}, \gamma_{23}$ be paths in $(0, \infty)$ 
from $z_{1}^{0}, z_{2}^{0}, z_{3}^{0}, z_{2}^{0}-z_{3}^{0}$, respectively, 
to $1$.
From the definitions of the module maps involved, 
we have the following equalities:
$$i'_{a_{2}}=\mathcal{T}_{\gamma_{1}}\circ i'_{a_{2}; z^{0}_{1}},$$
$$(I_{(W^{a_{2}})'}\boxtimes\; r^{-1}_{a_{2}})\circ 
\mathcal{T}_{\gamma_{1}}=
((I_{(W^{a_{2}})'}
\boxtimes \mathcal{T}_{\gamma_{2}})\circ \mathcal{T}_{\gamma_{1}})\circ
(I_{(W^{a_{2}})'}\boxtimes_{P(z^{0}_{1})}r^{-1}_{a_{2}, z^{0}_{2}})
$$
\begin{eqnarray*}
\lefteqn{(I_{(W^{a_{2}})'}\boxtimes (I_{W^{a_{2}}}
\boxtimes
i'_{a_{1}}))\circ ((I_{(W^{a_{2}})'}
\boxtimes \mathcal{T}_{\gamma_{2}})\circ \mathcal{T}_{\gamma_{1}})}\nn
&&= ((I_{(W^{a_{2}})'}
\boxtimes (I_{W^{a_{2}}} \boxtimes
\mathcal{T}_{\gamma_{3}}))\circ (I_{(W^{a_{2}})'}
\boxtimes \mathcal{T}_{\gamma_{2}})\circ \mathcal{T}_{\gamma_{1}})\circ\nn
&&\quad\quad\quad\quad\quad\quad\quad\quad\quad\quad
\circ (I_{(W^{a_{2}})'}\boxtimes_{P(z^{0}_{1})}(I_{W^{a_{2}}}
\boxtimes_{P(z^{0}_{2})}
i'_{a_{1}; z^{0}_{3}})),
\end{eqnarray*}
\begin{eqnarray*}
\lefteqn{(I_{(W^{a_{2}})'}\boxtimes
\mathcal{A})\circ ((I_{(W^{a_{2}})'}
\boxtimes (I_{W^{a_{2}}}\boxtimes 
\mathcal{T}_{\gamma_{3}}))\circ (I_{(W^{a_{2}})'}
\boxtimes \mathcal{T}_{\gamma_{2}})\circ \mathcal{T}_{\gamma_{1}})}\nn
&&=((I_{(W^{a_{2}})'}
\boxtimes (
\mathcal{T}_{\gamma_{23}}\boxtimes I_{W^{a_{2}}}))\circ (I_{(W^{a_{2}})'}
\boxtimes \mathcal{T}_{\gamma_{3}})\circ \mathcal{T}_{\gamma_{1}})\circ\nn
&&\quad\quad\quad\quad\quad\quad\quad\quad\quad\quad\quad\quad\quad\quad
\circ \left(I_{(W^{a_{2}})'}\boxtimes_{P(z^{0}_{1})} 
\mathcal{A}_{P(z^{0}_{2}), P(z^{0}_{2})}
^{P(z^{0}_{2}-z^{0}_{3}), P(z^{0}_{3})}\right),
\end{eqnarray*}
\begin{eqnarray*}
\lefteqn{(I_{(W^{a_{2}})'}\boxtimes
(\mathcal{C}^{2}\boxtimes
I_{W^{a_{2}}}))\circ (((I_{(W^{a_{2}})'}
\boxtimes 
(\mathcal{T}_{\gamma_{23}}\boxtimes I_{W^{a_{2}}}))\circ (I_{(W^{a_{2}})'}
\boxtimes \mathcal{T}_{\gamma_{3}})\circ \mathcal{T}_{\gamma_{1}})}\nn
&&=(((I_{(W^{a_{2}})'}
\boxtimes 
(\mathcal{T}_{\gamma_{23}}\boxtimes I_{W^{a_{2}}}))\circ (I_{(W^{a_{2}})'}
\boxtimes \mathcal{T}_{\gamma_{3}})\circ \mathcal{T}_{\gamma_{1}})\circ\nn
&&\quad\quad\quad\quad\quad\quad\quad\quad\quad\quad
\circ \left(I_{(W^{a_{2}})'}\boxtimes_{P(z^{0}_{1})} 
\left(\left(\mathcal{C}_{P(z^{0}_{2}-z^{0}_{3})}\right)^{2}
\boxtimes_{P(z^{0}_{3})} 
I_{W^{a_{2}}}\right)\right),
\end{eqnarray*}
\begin{eqnarray*}
\lefteqn{(I_{(W^{a_{2}})'}\boxtimes
\mathcal{A}^{-1})\circ (((I_{(W^{a_{2}})'}
\boxtimes 
(\mathcal{T}_{\gamma_{23}}\boxtimes I_{W^{a_{2}}}))\circ (I_{(W^{a_{2}})'}
\boxtimes \mathcal{T}_{\gamma_{3}})\circ \mathcal{T}_{\gamma_{1}})}\nn
&&=((I_{(W^{a_{2}})'}
\boxtimes (I_{W^{a_{2}}} \boxtimes
\mathcal{T}_{\gamma_{3}}))\circ (I_{(W^{a_{2}})'}
\boxtimes \mathcal{T}_{\gamma_{2}})\circ \mathcal{T}_{\gamma_{1}})\circ\nn
&&\quad\quad\quad\quad\quad\quad\quad\quad\quad\quad\quad
\circ \left(I_{(W^{a_{2}})'}\boxtimes_{P(z^{0}_{1})} 
\left(\mathcal{A}_{P(z^{0}_{2}), P(z^{0}_{2})}
^{P(z^{0}_{2}-z^{0}_{3}), P(z^{0}_{3})}\right)^{-1}\right),
\end{eqnarray*}
\begin{eqnarray*}
\lefteqn{(I_{(W^{a_{2}})'}\boxtimes
(I_{W^{a_{2}}}\boxtimes
\hat{e}_{a_{1}}))\circ
((I_{(W^{a_{2}})'}
\boxtimes (I_{W^{a_{2}}} \boxtimes
\mathcal{T}_{\gamma_{3}}))\circ (I_{(W^{a_{2}})'}
\boxtimes \mathcal{T}_{\gamma_{2}})\circ \mathcal{T}_{\gamma_{1}})}\nn
&&\quad\quad =((I_{(W^{a_{2}})'}
\boxtimes \mathcal{T}_{\gamma_{2}})\circ \mathcal{T}_{\gamma_{1}})
\circ (I_{(W^{a_{2}})'}\boxtimes_{P(z^{0}_{1})}
(I_{W^{a_{2}}}\boxtimes_{P(z^{0}_{2})}
\hat{e}_{a_{1}; z^{0}_{3}})),
\end{eqnarray*}
$$(I_{(W^{a_{2}})'}\boxtimes
r_{a_{2}})\circ ((I_{(W^{a_{2}})'}
\boxtimes \mathcal{T}_{\gamma_{2}})\circ \mathcal{T}_{\gamma_{1}})
= \mathcal{T}_{\gamma_{1}}\circ (I_{(W^{a_{2}})'}\boxtimes_{P(z^{0}_{1})}
r_{a_{2}; z^{0}_{2}}),$$
$$\hat{e}_{a_{2}; z^{0}_{1}}\circ \mathcal{T}_{\gamma_{1}}
=\hat{e}_{a_{2}}.$$
Each of these equalities is equivalent to a commutative diagram.
For example, the first and the second equalities above are
equivalent to
the commutative diagrams
$$\begin{CD}
V @>i'_{a_{2}; z^{0}_{1}}>>
(W^{a_{2}})'\boxtimes_{P(z^{0}_{1})} W^{a_{2}}\\
@V=VV @VV\mathcal{T}_{\gamma_{1}}V\\
V@>i'_{a_{2}}>> 
(W^{a_{2}})'\boxtimes W^{a_{2}}
\end{CD}
$$
and 
$$\begin{CD}
(W^{a_{2}})'\boxtimes_{P(z^{0}_{1})} W^{a_{2}}
@>I_{(W^{a_{2}})'}\boxtimes_{P(z^{0}_{1})}r^{-1}_{a_{2}, z^{0}_{2}}>>
(W^{a_{2}})'\boxtimes_{P(z^{0}_{1})} 
(W^{a_{2}}\boxtimes_{P(z^{0}_{2})} V)\\
@V\mathcal{T}_{\gamma_{1}}VV @VV(I_{(W^{a_{2}})'}
\boxtimes \mathcal{T}_{\gamma_{2}})\circ \mathcal{T}_{\gamma_{1}}V\\
(W^{a_{2}})'\boxtimes W^{a_{2}}
@>I_{(W^{a_{2}})'}\boxtimes\; r^{-1}_{a_{2}}>>
(W^{a_{2}})'\boxtimes
(W^{a_{2}}\boxtimes V)
\end{CD}
$$
(Since  the size of some of the diagram is too big, we leave the
other commutative
diagrams to the reader.)

Combining all the diagrams above or equivalently, all the
equalities above, we see that the composition of the module 
maps in the sequence (\ref{link}) and the composition of the module 
maps in the sequence (\ref{v-link}) are equal. Then by Proposition 
\ref{v-link-calc-1}, we obtain the conclusion of the 
proposition.
\epfv 

\begin{cor}
The composition of the module maps in the sequence obtained from 
(\ref{link}) by replacing $\hat{e}_{a_{1}}$ and $\hat{e}_{a_{1}}$ 
by $e_{a_{1}}$ and $e_{a_{1}}$, respectively,
is equal to 
$$\frac{(B^{(-1)})^{2}_{a_{2}, a_{1}}}{F_{a_{1}}F_{a_{2}}}I_{V},$$
that is, 
\begin{eqnarray}\label{modular-calc3}
\lefteqn{e_{a_{2}'}\circ 
(I_{(W^{a_{2}})'}\boxtimes
r_{a_{2}})
\circ (I_{(W^{a_{2}})'}\boxtimes
(I_{W^{a_{2}}}\boxtimes
e_{a'_{1}}))}\nn
&&\circ (I_{(W^{a_{2}})'}\boxtimes
\mathcal{A}^{-1})\circ (I_{(W^{a_{2}})'}\boxtimes
(\mathcal{C}^{2}
\boxtimes
I_{W^{a_{2}}}))
\circ (I_{(W^{a_{2}})'}\boxtimes
\mathcal{A})\nn
&&\circ
(I_{(W^{a_{2}})'}\boxtimes
(I_{W^{a_{2}}}\boxtimes
i_{a'_{1}}))
\circ (I_{(W^{a_{2}})'}\boxtimes
r^{-1}_{a_{2}})\circ i_{a'_{2}}
\nn
&&=\frac{(B^{(-1)})_{a_{2}, a_{1}}^{2}}
{F_{a_{1}}
F_{a_{2}}}I_{V}.
\end{eqnarray}
\end{cor}
\pf
From (\ref{modular-calc2}) and the definitions of 
$e_{a_{1}}$ and $e_{a_{2}}$, 
we obtain (\ref{modular-calc3}).
\epfv

\begin{thm}\label{nondeg}
Let $V$ be a simple vertex operator algebra satisfying the 
conditions in Section 1.
Then the ribbon category structure on the 
category of $V$-modules constructed in 
\cite{HL1}--\cite{HL4}, \cite{H1} and \cite{H5}
is nondegenerate.
\end{thm}
\pf
By definition, the traces of $\mathcal{C}^{2}\in \hom((W^{a_{1}})', 
W^{a_{2}})$ for $a_{1}, a_{2}\in \mathcal{A}$ can be 
calculated as follows: For any $a_{1}, a_{2}\in \mathcal{A}$,
consider the composition of the module maps 
in the sequence obtained from (\ref{link}) by replacing 
$\hat{e}_{a_{1}}$ and $\hat{e}_{a_{1}}$ 
by $e_{a_{1}}$ and $e_{a_{1}}$, respectively, that is,
consider the module map given by the left-hand side of 
(\ref{modular-calc3}). This is a module map from $V$ to $V$. 
Since $V$ is an irreducible $V$-module, the module map we 
are considering is equal to the identity operator on $V$ multiplied by 
a number. This number is the trace of $\mathcal{C}^{2}\in \hom((W^{a_{1}})', 
W^{a_{2}})$.

Thus by Proposition \ref{v-link-calc-2}, 
the traces of $\mathcal{C}^{2}\in \hom((W^{a_{1}})', 
W^{a_{2}})$ for $a_{1}, a_{2}\in \mathcal{A}$
are
$$\frac{(B^{(-1)})_{a_{2}, a_{1}}^{2}}
{F_{a_{1}}
F_{a_{2}}}.
$$
By (\ref{s-form-3}), they are actually equal to 
$$\frac{S_{a_{1}}^{a_{2}}}{S_{e}^{e}}, \;a_{1}, a_{2}\in \mathcal{A},$$
which form an invertible matrix.  So the 
tensor category is nondegenerate.
\epfv

From Theorem \ref{nondeg} and the definition of 
modular tensor category (see for example \cite{T} and 
\cite{BK}), we immediately obtain the main result of the 
present paper:

\begin{thm}
Let $V$ be a simple vertex operator algebra satisfying the 
conditions in Section 1.
Then the category of $V$-modules has a natural structure
of modular tensor category.
\end{thm}

\begin{rema}
{\rm As in the discussion in Remark \ref{v-rigidity}, 
we can also introduce a notion of nondegeneracy for semisimple 
rigid vertex tensor categories and a notion of modular 
vertex tensor category. In fact, for any semisimple rigid
vertex tensor category, the sequence (\ref{v-link}) still makes sense. So
we can introduce the notion of nondegeneracy for vertex tensor categories
by requiring the compositions of the morphisms in this sequence
to be equal to the identity on the unit object. Twists and balancing axioms
for vertex tensor categories can also be introduced in an 
obvious way. A 
modular vertex tensor category is a semisimple vertex tensor 
category which is also rigid, balanced and nondegenerate.
Then our result above shows that 
the vertex tensor category of $V$-modules is nondegenerate 
and thus the category has a natural structure of modular vertex tensor 
categories. Moreover, the proof of Proposition \ref{v-link-calc-2}
actually shows that if a  semisimple rigid vertex tensor 
category is nondegenerate, then the corresponding semisimple
rigid braided tensor category is also nondegenerate. In particular,
a modular vertex tensor category gives a modular tensor category.}
\end{rema}

\noindent {\small \sc Department of Mathematics, Rutgers University,
110 Frelinghuysen Rd., Piscataway, NJ 08854-8019}

\noindent {\em E-mail address}: {\tt yzhuang@math.rutgers.edu}

\end{document}